\documentclass[12pt]{article}

\usepackage[utf8]{inputenc}
\usepackage[T1]{fontenc}

\usepackage{geometry}
\usepackage{ragged2e}
\usepackage[pagebackref]{hyperref}
\hypersetup{
colorlinks=true, 
citecolor = blue, 
linkcolor = blue, 
urlcolor = blue, 
anchorcolor = blue
}

\usepackage{amsmath, amsthm, amsfonts}
\usepackage{amscd, amssymb, latexsym}
\usepackage{soul}
\usepackage{mathabx, upgreek}

\usepackage[all]{xy}
\usepackage{mathtools} 
\usepackage{caption}

\usepackage{comment, color}
\usepackage{enumitem}

\usepackage{paratext/commands_Hexp}

\usepackage{mathrsfs}
\usepackage{array}
\usepackage{diagbox}
\usepackage{graphicx}

\usepackage{tikz, pgf, pgfplots}
\pgfplotsset{compat=newest}

\usetikzlibrary[patterns]
\usetikzlibrary{arrows, arrows.meta}
\usetikzlibrary{calc} 
\def\centerarc[#1](#2)(#3:#4:#5)%
    {\draw [#1] ($(#2)+({#5*cos(#3)}, {#5*sin(#3)})$) arc (#3:#4:#5);} %

\usetikzlibrary{intersections,through,backgrounds} 
\usetikzlibrary{shapes.geometric} 

\definecolor{grey}{rgb}{0.5,0.5,0.5}
\definecolor{forestgreen}{rgb}{0.133,0.545,0.133}
\definecolor{marron}{rgb}{0.6,0.2,0.}
\definecolor{violet}{rgb}{0.5,0.,0.5}
\definecolor{lightpurple}{rgb}{.6,.2,.6}

\usepackage[outline]{contour}

\usepackage{subfiles}

\title{Geometry and Transcendence of the Hexponential}

\author{Scott Schmieding, Christopher-Lloyd Simon}
\date{\today}

\begin{document}

\maketitle

\begin{abstract}
The modular group $\PSL_2(\Z)$ acts on the upper-half plane $\HP$ with quotient the modular orbifold, uniformized by the $\J$-function $\J\colon \HP\to \C$.
%
We first show that second derived subgroup $\PSL_2(\Z)''$ corresponds to a $\Z^2\rtimes \Z/6$ Galois cover of the modular orbifold by a hexpunctured plane, uniformized by the hexponential map $\hexp \colon \HP \to \C\setminus (\omega_0\Z[j])$, which is a primitive of $C\eta^4$ where $\omega_0\in i\R$ and $C\in \R$ are explicit constants and $\eta$ is the Dedekind eta function.
We describe the values of the cusp-compactification $\partial \hexp\colon \Q\P^1\to \omega_0 \Z[j]$.
After defining the radial-compactification $\Shexp \colon \Radial \to \R/(2\pi\Z)$, we construct a simple section $\InSh \colon \R/(2\pi\Z) \to \Sturm \bmod{\PSL_2(\Z)'}$ where $\Sturm \subset \R\P^1$ is a set of numbers whose continued fraction expansions arise from Sturmian sequences, which contains the set $\Markov$ of Markov quadratic irrationals as those numbers arising from periodic Sturmian sequences. 
We will show that the values of $\InSh$ are either Markov quadratic irrationals or transcendental.
Finally we provide a continued fraction expansion for $\hexp$, and discuss its monodromy.
\end{abstract}

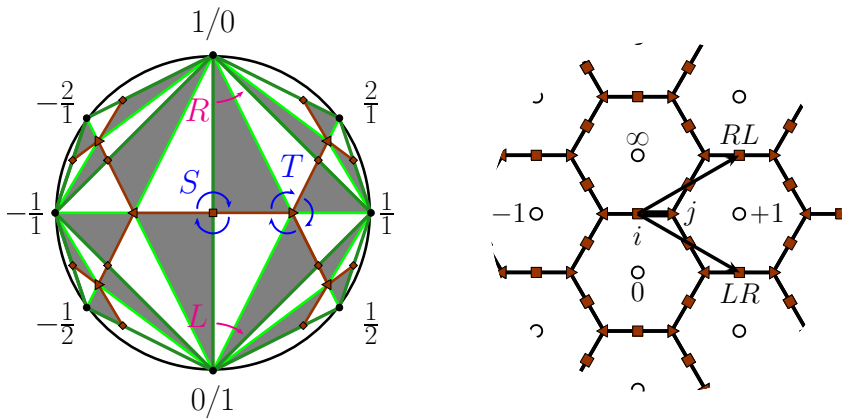
\begin{figure}[h]
    \centering
    \scalebox{0.5}{\subfile{images/tikz/PSL2Z-pavage-HP_1}}
    \quad
    \scalebox{0.9}{\begin{tikzpicture}[rotate=30]
\clip[shift={(2.5,2.5)}] (-30:1) circle (26mm);

\begin{scope}[shift={(2*2*0.8657,2*2*0.75)}]
    \draw[color=black, line width=3pt] (-120:0.8657) -- (-90:1);
    \draw[color=black, line width=3pt, anchor=north] (-120:0.8657) node {$i$};
    \draw[color=black, line width=3pt, anchor=west] (-90:1) node {$j$};
\end{scope}

\foreach \i in {0,...,3}
\foreach \j in {1,3} {
\begin{scope}[shift={(2*\i*0.8657+0.8657,2*\j*0.75)}]
\draw [line width=1.5pt,color=black] plot[samples at={-150,-90,...,210},variable=\x] (\x:1);
\draw [fill=marron, rotate around={45:(0:0.8657)}] (0:0.8657) ++(-3.pt,0 pt) -- ++(3.pt,3.pt)--++(3.pt,-3.pt)--++(-3.pt,-3.pt)--++(-3.pt,3.pt);
\draw [fill=marron, rotate around={45:(-180:0.8657)}] (-180:0.8657) ++(-3.pt,0 pt) -- ++(3.pt,3.pt)--++(3.pt,-3.pt)--++(-3.pt,-3.pt)--++(-3.pt,3.pt);
\draw [line width=0.8pt,color=black] (0,0) circle (2.5pt);
\end{scope}
}

\foreach \i in {0,...,3}
\foreach \j in {0,2} {
    \begin{scope}[shift={(2*\i*0.8657,2*\j*0.75)}]
        \draw [line width=1.5pt,color=black] plot[samples at={-150,-90,...,210},variable=\x] (\x:1);
        \draw [fill=marron, rotate around={45:(-180:0.8657)}] (-180:0.8657) ++(-3.pt,0 pt) -- ++(3.pt,3.pt)--++(3.pt,-3.pt)--++(-3.pt,-3.pt)--++(-3.pt,3.pt);
        \draw [fill=marron, rotate around={15:(-120:0.8657)}] (-120:0.8657) ++(-3.pt,0 pt) -- ++(3.pt,3.pt)--++(3.pt,-3.pt)--++(-3.pt,-3.pt)--++(-3.pt,3.pt);
        \draw [fill=marron, rotate around={-15:(-60:0.8657)}] (-60:0.8657) ++(-3.pt,0 pt) -- ++(3.pt,3.pt)--++(3.pt,-3.pt)--++(-3.pt,-3.pt)--++(-3.pt,3.pt);
        \draw [fill=marron, rotate around={45:(0:0.8657)}] (0:0.8657) ++(-3.pt,0 pt) -- ++(3.pt,3.pt)--++(3.pt,-3.pt)--++(-3.pt,-3.pt)--++(-3.pt,3.pt);
        \draw [fill=marron, rotate around={15:(60:0.8657)}] (60:0.8657) ++(-3.pt,0 pt) -- ++(3.pt,3.pt)--++(3.pt,-3.pt)--++(-3.pt,-3.pt)--++(-3.pt,3.pt);
        \draw [fill=marron, rotate around={-15:(120:0.8657)}] (120:0.8657) ++(-3.pt,0 pt) -- ++(3.pt,3.pt)--++(3.pt,-3.pt)--++(-3.pt,-3.pt)--++(-3.pt,3.pt);
        \draw [fill=marron,shift={(-150:1)},rotate=180] (0,0) ++(0 pt,{4.5/1.5pt}) -- ++({3.897/1.5pt},{-6.75/1.5pt})--++({-7.794/1.5pt},0 pt) -- ++({3.897/1.5pt},{6.75/1.5pt});
        \draw [fill=marron,shift={(-90:1)},rotate=0] (0,0) ++(0 pt,{4.5/1.5pt}) -- ++({3.897/1.5pt},{-6.75/1.5pt})--++({-7.794/1.5pt},0 pt) -- ++({3.897/1.5pt},{6.75/1.5pt});
        \draw [fill=marron,shift={(-30:1)},rotate=180] (0,0) ++(0 pt,{4.5/1.5pt}) -- ++({3.897/1.5pt},{-6.75/1.5pt})--++({-7.794/1.5pt},0 pt) -- ++({3.897/1.5pt},{6.75/1.5pt});
        \draw [fill=marron,shift={(30:1)},rotate=0] (0,0) ++(0 pt,{4.5/1.5pt}) -- ++({3.897/1.5pt},{-6.75/1.5pt})--++({-7.794/1.5pt},0 pt) -- ++({3.897/1.5pt},{6.75/1.5pt});
        \draw [fill=marron,shift={(90:1)},rotate=180] (0,0) ++(0 pt,{4.5/1.5pt}) -- ++({3.897/1.5pt},{-6.75/1.5pt})--++({-7.794/1.5pt},0 pt) -- ++({3.897/1.5pt},{6.75/1.5pt});
        \draw [fill=marron,shift={(150:1)},rotate=0] (0,0) ++(0 pt,{4.5/1.5pt}) -- ++({3.897/1.5pt},{-6.75/1.5pt})--++({-7.794/1.5pt},0 pt) -- ++({3.897/1.5pt},{6.75/1.5pt});
        \draw [line width=0.8pt,color=black] (0,0) circle (2.5pt);
    \end{scope}
    }
    

\draw[color=black,anchor=east] (2*1*0.8657,2*2*0.75) node{$-1$};
\draw[color=black,anchor=west] (2*2.5*0.8657,2*1*0.75) node{$+1$};
\draw[color=black,anchor=south] (2*2*0.8657,2*2*0.75) node {$\infty$};
\draw[color=black,anchor=north] (2*1*0.8657+0.8657,2*1*0.75) node {$0$};

\begin{scope}[shift={(2*2*0.8657,2*2*0.75)}]
    \draw [color=black,line width=1.5pt, -{Stealth[length=2.mm,width=1.5mm]}, anchor=south] (-120:0.8657) -- (-30:2*0.75) node {$RL$};
    \draw [color=black,line width=1.5pt, -{Stealth[length=2.mm,width=1.5mm]}, anchor = north] (-120:0.8657) -- (-79:2.3) node{$LR$};
\end{scope}
\end{tikzpicture}}
\end{figure}


\newpage

\renewcommand{\contentsname}{Plan of the paper}
\setcounter{tocdepth}{2}
\tableofcontents

\section{Introduction}

The paper unfolds the plan of the introduction. 

\subsection{The modular orbifold and euclidean continued fractions}

The modular group $\Gamma=\PSL_2(\Z)$ acts by conformal transformations of the hyperbolic plane $\HP=\{z\in \C\mid \Im(z)>0\}$ with boundary $\partial \HP = \R\P^1$, as the stabiliser of its ideal triangulation $\triangle$ with vertices $\Q\P^1$ and edges $\left(\frac{a}{c},\frac{b}{d}\right)$ such that $ad-bc=1$. Let
\begin{equation*}\textstyle
    S=\begin{psmallmatrix}
    0 & -1\\
    1 & 0
    \end{psmallmatrix}
    \quad
    T =
    \begin{psmallmatrix}
    1 & -1\\
    1 & 0
    \end{psmallmatrix}
    \qquad
    L = T^{-1}S = 
    \begin{psmallmatrix}
    1 & 0\\
    1 & 1
    \end{psmallmatrix}
    \quad
    R = TS^{-1} =
    \begin{psmallmatrix}
    1 & 1\\
    0 & 1
    \end{psmallmatrix}.
\end{equation*}

The modular orbifold $\M = \Gamma \backslash \HP$ has genus zero, a cusp associated to the fixed point $\infty \in \Q\P^1 \subset \partial \HP$ of $R$, as well as two conical singularities of order $2$ and $3$ associated to the fixed points $i,j \in \HP$ of $S$ and $T$.
Thus $\Gamma = \pi_1(\M)$ is the free amalgam of its subgroups $\Z/2$ and $\Z/3$ generated by $S$ and $T$. 

The elements of $\Gamma$ stabilising $[0,+\infty]\subset \R\P^1$ form the euclidean monoid $\PSL_2(\N)$, which is freely generated by the parabolic matrices $R$ and $L$.
The subset of sequences in $\{L,R\}^\N$ which do not end with an infinite string of $L$ are in bijective correspondence with the continued fraction expansions of positive numbers in $(0,+\infty]\subset \R\P^1$. 
The following result can be traced back to the works of Gauss \cite{Gauss_disquisitiones_1807} and Galois \cite{Galois_fraction-continue_1828}.

\begin{Proposition}[Geodesics in $\M$] 
\label{prop:Intro_Gauss-Galois-contfrac}
A geodesic in $\M$ has a lift $(\alpha_-,\alpha_+)\subset \HP$ whose endpoints satisfy $-1<\alpha_-<0$ and $1<\alpha_+<\infty$.
It intersects $\triangle$ along a sequence of triangles whose $\{L,R\}^{\Z}$-encoding is obtained from the continued fraction expansions of $\alpha_+$ and $-1/\alpha_-$ by concatenating the transpose of the latter with the former.

The closed geodesics in $\M$ correspond to the periodic sequences, hence to the $\Gamma$-orbits of pairs $(\alpha_-,\alpha_+)$ of Galois conjugate quadratic irrationals.
\end{Proposition}

\subsection{The modular torus and Sturmian numbers}

The derived subgroup $\Gamma'$ of $\Gamma$, which is freely generated by $X = LR$ and $Y = RL$, corresponds to a Galois cover of $\M$ by a punctured torus $\M'$ with group $\Gamma/\Gamma'=\Z/6$.
The simple geodesics of $\M'$ were described by Cohn \cite{Cohn_Markoff-perforated-torus_1971} and Series \cite{Series_Geo-Markov-Num_1985}.

\begin{Proposition}[Simple geodesics in $\M'$]
\label{prop:Intro_Markov-Series-contfrac}
    Fix $-1/\alpha_-, \alpha_+\in \R_+^*$ and consider the $L\& R$-sequence obtained by concatenating their (transposed) continued fractions.
    
    The geodesic $(\alpha_-,\alpha_+)\subset \HP$ projects to a geodesic in $\M'$ which is simple if and only if the $L\& R$ can be paired to form an $X\& Y$-sequence which is Sturmian.
    
    The simple closed geodesics in $\M'$ correspond combinatorially to the periodic Sturmian $X\& Y$-sequences, and arithmetically to the Markov quadratic irrationals.
\end{Proposition}

The conformal structure on $\M'$ is given by a punctured elliptic curve isomorphic over $\R$ to the affine cubic $y^2=(x-\omega_0 j^0)(x-\omega_0j^1)(x-\omega_0 j^2)$ with abelian differential $dx/y=-dx/\sqrt{x^3-\rvert \omega_0\rvert ^3}$ having fundamental period $\omega_0=\frac{2\pi^{1/2}}{i3^{1/4}}$.

This abelian differential on $\M'$ lifts on the universal cover $\HP$ to $C\eta^4(\tau)d\tau$ for some $C\in \C^*$ where $\eta$ is the Dedekind eta function, whose expansion in $q=e^{i2\pi \tau}$:
    \begin{equation*}
        \eta(\tau)=\sum_{n=0}^\infty \chi(n)q^{\tfrac{1}{24}n^2}
        \quad \mathrm{yields} \quad
        \eta^4(\tau)=\sum_{n=0}^\infty \psi(n)q^{\tfrac{1}{24}n}
    \end{equation*}
where $\chi$ is the unique primitive character $\bmod{12}$ (that is the Jacobi symbol $\bmod{12}$, equal to $\chi(12n\pm 1)=1$, $\chi(12n\pm 5)=-1$ and $\chi(n)=0$ if $2\mid n$ or $3\mid n$), and where $\psi(n)$ is the sum of $\chi(abcd)$ over $\{(a,b,c,d)\in \N^4 \mid a^2+b^2+c^2+d^2=n\}$.

\subsection{The hexpunctured plane and its cusp compactifications}

The second derived subgroup $\Gamma''\subset \Gamma'$ corresponds to the universal abelian cover $\M''\to \M'$ of the punctured torus by a lattice-punctured plane with Galois group $\Gamma'/\Gamma''=\Z_X\oplus \Z_Y$.
The action of the Galois group $\Gamma/\Gamma''$ on the corresponding cover $\M''\to \M$ was described in the second author's thesis \texorpdfstring{\cite[Section 3.2]{CLS_phdthesis_2022}}{CLSthesis}.

\begin{Proposition}[Hexagonal lattice]
\label{prop:Intro_CLS_hexagonal-group}
We have a semi-direct decomposition 
\[\Gamma/\Gamma''=\Gamma'/\Gamma'' \rtimes \Gamma/\Gamma'\] 
where the action of the quotient $\Gamma/\Gamma'=\Z/2\times \Z/3$ by outer-automophisms of the kernel $\Gamma'/\Gamma''=\Z_X\oplus \Z_Y$ is given by $\AD_S=D_{S^2}=D_{-\Id}$ and $\AD_T=D_{T^2}$ since we have $SXS^{-1}=X^{-1}$, $SYS^{-1}=Y^{-1}$ and $TXT^{-1}=Y^{-1}$, $TYT^{-1}=-XY^{-1}$. 

This represents $\Gamma/\Gamma''$ as the affine isometry group of the oriented hexagonal lattice $H_1(\Gamma';\Z)=\Z_X\oplus \Z_Y$ with $angle{(X,Y)}=\frac{2\pi}{6}$, where the translation action of $\Gamma'/\Gamma''$ is by $\Z^2$-translation while the outer-automorphism action of $\Gamma/\Gamma'$ is by $\Z/6$-rotation.    
\end{Proposition}

The integration of abelian differential on $\M'$ identifies its universal abelian cover $\M''$ with $H_1(\M';\R)\setminus H_1(\M';\Z)$, hence with the hexpunctured plane $\C\setminus \omega_0\Z[j]$.

The map $\hexp \colon \HP \to \M''$ in the coordinates $\{z\in \C \mid \Im(z)>0\} \to \C\setminus \omega_0 \Z[j]$ is therefore the primitive of $C\eta(z)^4dz$, with expansion:
\begin{equation*}
    \hexp(\tau)=\int_\infty^\tau C\eta^4(z)dz = \frac{12C}{i\pi}\sum_{n=1}^\infty \frac{\psi(n)}{n}\cdot q^{\tfrac{1}{24}n}
\end{equation*}
We deduce the value of $C = \frac{2^{10/3}}{3^{3/4}} \frac{\pi^{5/2}}{\Upgamma(1/3)^{3}}$ from the relation $\omega_0= \int_\infty^0 C\eta^4(\tau)d\tau$.

The main results in this paper concern arithmetic and geometric properties of the uniformizing map $\hexp\colon \HP\to \M''$, and its compactifications.

We first describe the values of the cusp compactification $\partial \hexp \colon \Q\P^1\to \omega_0\Z[j]$.

\begin{Theorem}[Values of \texorpdfstring{$\partial \hexp$}{dhexp}]
\label{thm:Intro_values-hexp(rational)}
    Consider $a/c\in \Q\P^1$ with $\gcd(a,b)=1$.
    
    There exists $A=\begin{psmallmatrix} a & b \\ c & d \end{psmallmatrix} \in \Gamma'$, and we may consider its class $A\equiv X^mY^n \mod{\Gamma''}$. 
    The value $\partial \hexp(\frac{a}{c})$ of the improper integral and conditionally convergent series:
    \begin{equation*}
        \partial \hexp(\tfrac{a}{c}) = 
        \int_\infty^{\frac{a}{c}} C\eta^4(z)dz 
        = \tfrac{12C}{i\pi}\sum_{n=1}^\infty \tfrac{\psi(n)}{n}\cdot \exp\left(\tfrac{i\pi}{12}\tfrac{a}{c}n\right)
        \quad \text{where} \quad
        C = \tfrac{2^{10/3}}{3^{3/4}} \tfrac{\pi^{5/2}}{ \Upgamma(1/3)^{3}} 
    \end{equation*}
    \begin{equation*}
        \text{is equal to}\quad 
        \partial \hexp(\tfrac{a}{c}) = 
        \lvert \omega_0 \rvert \left(m\exp(-\tfrac{i\pi}{6})+n\exp(+\tfrac{i\pi}{6})\right)
        \quad \text{where} \quad
        \lvert \omega_0 \rvert = \tfrac{2\pi^{1/2}}{3^{1/4}}.
    \end{equation*}
\end{Theorem}

The function $\partial \hexp \colon \Q\P^1\to \C$ is a modular symbol for the group $\Gamma''$. Modular symbols were introduced by Manin in \cite{Manin_parabolic-points-modular-curves_1972} for the groups $\Gamma_0(N)$, and have since been the source of various generalisations and extensive studies (see \cite{Manin-Marcolli_Contfrac-modular-symbols-noncommutative-geometry_2002, Manin_lectures-modular-symbols_2009}).

\begin{Remark}[Rational correlations]
    The conditional convergence of the series $\partial\hexp(r)$ at rationals $r\in \Q$ implies a strong correlation between $n\mapsto \psi(n)$ and periodic rotations of the circle $n\mapsto n+r \bmod{1}$. 
    These correlations factor through the cosets $r\in \Gamma''\backslash \Q\P^1$, and their values are organised according to the lattice $\omega_0\Z[j]$.

    This may be equivalent to some kind of "Dirichlet equidistribution theorem" for the arithmetic progressions represented by a $\chi$-twisted version of the quadratic form $a^2+b^2+c^2+d^2$ or the divisor function $\sum_{d\mid n} d$.
\end{Remark}

\subsection{Transcendent arguments of the hexponential map}

We then define the ray compactification $\Shexp \colon \Radial \to \S H_1(\M';\R)$. 
Let $\Radial$ be the set of numbers $\alpha \in \R\P^1$ such that as $\tau \in \HP$ diverges to $\alpha \in \partial \HP$, the modulus $\lvert \hexp(\tau) \rvert \in \R_+$ diverges while the argument $\arg \hexp(\tau) \in \R/(2\pi \Z)$ converges.
Thus $\alpha \in \Radial$ when the geodesic $(i,\alpha) \subset \HP$ maps under $\hexp$ to a geodesic of $\M''$ which escapes towards a definite direction in $\S H_1(\M';\R)$: that is $\Shexp(\alpha)$.

\begin{Remark}[Singularities along the boundary]
We have $\Radial \subset \R \setminus \Q$.
If $\alpha \in \R\P^1$ is neither rational nor quadratic then $\hexp(\tau)$ diverges as $\tau\in \HP$ diverges to $\alpha\in \R\P^1$.

Suppose $\alpha \in \R\P^1$ is a fixed point of $A\in \Gamma$: if $A^6\notin \Gamma''$ then $\alpha \in \Radial$, but if $A^6\in \Gamma''$ then $\alpha \notin \Radial$ yet $\lim_\alpha \hexp(\tau) = \partial\hexp(r)$ for some $r\in \Q\P^1$.

The subset $\Radial \subset \R\P^1$ is $\Gamma'$-invariant, $\Gamma'/\Gamma$ equivariant, and dense.

The map $\Shexp \colon \Radial \to \S H_1(\M',\R)$ is uniformly continuous and surjective.
It sends sends uncountably many distinct $\Gamma'$-orbits in $\Radial$ to the same slope. 
\end{Remark}

The surjective maps $\Shexp \colon \Radial \to \S H_1(\M';\R)$ and $\partial \hexp \colon \Q\P^1 \to H_1(\M';\Z)$ have disjoint domains. Together, they define (for obvious topologies) a continuous compactification of $\hexp \colon \HP \to \M''$, equivariant under the left actions of $\Gamma$ and $\Gamma/\Gamma''$.
The translation action of $\Gamma'$ on $\S H_1(\M';\R)$ is trivial, so we have an equivariant map $\Shexp \colon \Radial \bmod{\Gamma'} \to \S H_1(\M';\R)$ under the left actions of $\Gamma$ and $\Gamma/\Gamma'$.

Finally, we define a simple section $\InSh \colon \S H_1(\M';\R) \to \Sturm \bmod{\Gamma'}$ of $\Shexp$ which is $\Gamma/\Gamma'$-equivariant and injective. A slope $(x,y) \in \S H_1(\M';\R)$ corresponds to a complete simple geodesic $\alpha\subset \M'$: it lifts in $\HP$ to the $\Gamma'$-orbit of a geodesic $(\alpha_-,\alpha_+)$ and we let $\InSh(x,y)=\Gamma'\cdot \alpha_+$.

We denote by $\Sturm$ the set of \emph{Sturmian} numbers $-1/\alpha_-,\alpha_+$ obtained in this way (namely whose continued fractions arise from a Sturmian $X^{\pm 1}\& Y^{\pm 1}$-sequence by replacing $X=LR$ and $Y=RL$).
In contains the subset $\Markov\subset \Sturm$ of \emph{Markov} numbers as those arising from rational slopes $(x,y) \in \S H_1(\M';\Z)$, namely those with periodic $X^{\pm 1}\& Y^{\pm 1}$ sequences (coinciding with the Markov quadratic irrationals).

\begin{Remark}[Cantor set]
    The set $\Sturm \subset \R\P^1$ is compact, with a dense countably infinite set of isolated points whose complement is a Cantor set. By \cite{Pelczynski_zero-dimensional-spaces_1965} such properties describe a unique space up to homeomorphism.
    It has Hausdorff dimension zero.
    The map $\Shexp \colon \Sturm \to \S H_1(\M';\R)$ is continuous and surjective.
\end{Remark}

\begin{Question}[Sturmian correlations]
    We observe that the sequence $\psi(n)/n$ bears very strong (lattice-like) correlations with rational numbers (bounded complexity), and rather strong (Cantor-like) correlations with Sturmian numbers (complexity $l+1$), and some correlations with the numbers in $\Radial$.

    We expect that one may characterise the subset of Sturmian numbers $\Sturm \subset \Radial$ by quantifying how fast $\arg \hexp(\theta)$ converges to $\Shexp(\alpha)$.
    We wonder if one can provide any arithmetic or diophantine interpretation of the numbers in $\Radial$.
\end{Question}


\begin{Theorem}[Transcendence of the Inshection map]
\label{thm:Intro_transcendence-hexp}
    For $(x,y) \in \S H_1(\M;\R)$:
    \begin{itemize}[noitemsep]
        \item If $(x,y) \in \S H_1(\M';\Z)$ is rational then $\InSh(x,y)\subset \Markov$ are quadratic.
        \item If $(x,y) \notin \S H_1(\M';\Z)$ is irrational then $\InSh(x,y)\subset \Sturm \setminus \Markov$ are transcendent.
    \end{itemize}
\end{Theorem}

\begin{Question}
Can one find continued fraction expansions for the real functions $\tan \circ \Shexp \colon \Radial \to \R\R^1$ and $\InSh \circ \arctan \colon \R\P^1 \to \Sturm$?
\end{Question}

We will however exhibit a continued fraction expansion for $\hexp$ in terms of the modular function $\lambda \colon \HP \to \C\setminus \{0,1,\infty\}$ (uniformizing the congruence cover $\M(2)$ of $\M$ associated to the congruence subgroup $\Gamma(2)$ of $\Gamma$, kernel of $\PSL_2(\Z) \to \PSL_2(\Z/2)$), which was kindly communicated to us by Bill Duke.

The Section \ref{sec:other-surfaces} announces results and proposes conjectures about other surfaces.
The appendix \ref{sec:more-on-psi} records several well known arithmetic properties of $\psi$.

\newpage

\section{The modular orbifold and continued fractions}

\subsection{Modular group and its conjugacy classes}

The modular group $\Gamma = \PSL_2(\Z)$ acts on the ideal triangulation $\triangle$ of $\HP$ with vertex set $\Q\P^1$ and edges $\left(\frac{a}{c},\frac{b}{d}\right)$ such that $ad-bc=1$, freely transitively on its flags, corresponding to the half-edges of its dual trivalent tree $\TT$.
The connected components of $\HP \setminus \TT$ correspond to the vertices of $\triangle$, parametrized by $\Gamma/ \langle R \rangle$.

\begin{equation*}
    S =
    \begin{pmatrix}
    0 & -1\\
    1 & 0
    \end{pmatrix}
    \quad
    T =
    \begin{pmatrix}
    1 & -1\\
    1 & 0
    \end{pmatrix}
    \qquad
    L = 
    \begin{pmatrix}
    1 & 0\\
    1 & 1
    \end{pmatrix}
    \quad
    R = 
    \begin{pmatrix}
    1 & 1\\
    0 & 1
    \end{pmatrix}
\end{equation*}
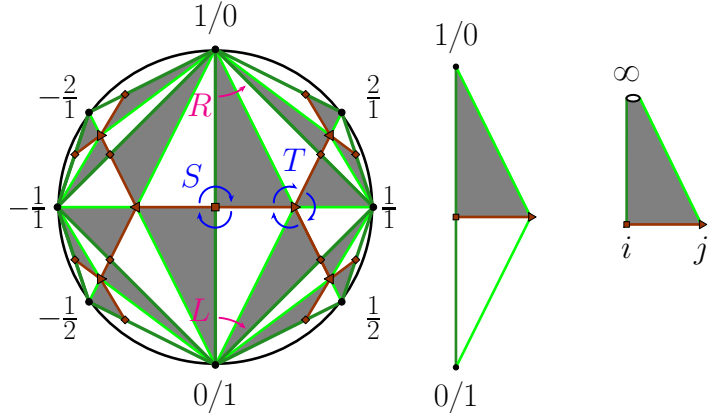
\begin{figure}[h]
    \centering
    \scalebox{0.5}{\subfile{images/tikz/PSL2Z-pavage-PH}}
    %
    %
\end{figure}

The modular orbifold $\M = \Gamma \backslash \HP$ has genus zero, a cusp associated to the fixed point $\infty \in \Q\P^1 \subset \partial \HP$ of $R$, as well as two conical singularities of order $2$ and $3$ associated to the fixed points $i,j \in \HP$ of $S$ and $T$.

Thus $\Gamma = \pi_1(\M)$ is the free amalgam of its subgroups $\Z/2$ and $\Z/3$ generated by $S$ and $T$. 
(Hence $\SL_2(\Z)$ is the amalgam of its subgroups $\Z/4$ and $\Z/6$ generated by $S$ and $T$ over their intersection $\Z/2$ generated by $S^2=-\Id=T^3$.)

It follows that in $\Gamma$, a finite order element is conjugate to a power of $S$ or $T$; whereas an infinite order $A$ is conjugate to $\prod_{1}^{l} T^{\epsilon_i}S^{-\epsilon_i}$ for a unique $l=\len A\in \N^*$ and a unique sequence of $\epsilon_i = \pm 1$ up to cyclic permutation: this corresponds to cyclic word in $R=TS^{-1}$ and $L=T^{-1}S$, and its Rademacher number is defined by 
\begin{equation*}
    \Rad(A)= \sum_{i=1}^{\len(A)} \epsilon_i = \#R-\# L
\end{equation*}


Conjugacy classes in $\Gamma$ correspond to orbifold homotopy classes of loops in $\M$.
The elliptic and parabolic classes correspond to loops circling around the singularities and cusp, whereas every hyperbolic class is represented by a unique closed geodesic.

An element of $\Gamma$ is primitive when it generates a maximal cyclic subgroup. This notion is invariant by conjugacy. Correspondingly, a cyclic $L\& R$-word is primitive when it is not a proper power, and a closed geodesic is primitive when it does not wind several times around itself.

A complete geodesic of $\HP$ is uniquely determined by the sequence of triangles in $\triangle$ that it intersects: it corresponds to a geodesic of the dual tree $\TT$ which is not horocyclic (bounding a region of the complement $\HP\setminus \TT$).
Such geometric and combinatorial geodesics are equally determined by their common endpoints in $\R\P^1$.

A closed geodesic of $\M$ lifts to the geodesics in $\HP$ corresponding to the periodic geodesics of $\TT$ whose period is given by $L\& R$-cycle of the associated conjugacy class.
Their endpoints form the $\Gamma$-orbit of a pair of Galois-conjugate quadratic irrationals.
Such geodesics determine a unique hyperbolic conjugacy class in $\Gamma$ which is primitive. 

\subsection{Euclidean monoid and continued fraction expansions}

The monoid $\SL_2(\N)\subset \SL_2(\Z)$ is freely generated by the transvections $L$ and $R$, and
it identifies with its image $\PSL_2(\N)\subset \PSL_2(\Z)$ which we call the euclidean monoid.

In $\PSL_2(\Z)$, the conjugacy classes of torsion elements are those of $\Id,S,T,T^{-1}$, whereas the conjugacy class of an infinite order element intersects the euclidean monoid $\PSL_2(\N)$ along the cyclic permutations of a unique $L\&R$-word.
Hence the hyperbolic geodesics of $\M$ are indexed by cyclic words in $L\& R$ containing both letters.

Observe that for $A\in \SL_2$, its transpose $A^\dag$ and inverse $A^{-1}$ are conjugate by $S$. Hence the transposition of cyclic words in $\PSL_2(\N)$ corresponds to inversion of the associated conjugacy classes in $\PSL_2(\Z)$.


Every positive real number $x$ admits a Euclidean continued fraction expansion:
\begin{equation*}
\Ecf{n_0,n_1,\dots}=n_0+\frac{1}{n_1+\dots}
\quad \mathrm{with} 
\quad n_j \in \N
\quad \mathrm{and} 
\quad \forall j>0,\; n_j>0.
\end{equation*}
Such an expansion is infinite if and only if $x$ is irrational, in which case it is unique. A rational $x$ has two expansions $x=\Ecf{n_0,\dots,n_k}$: one for which $n_k = 1$, the other for which $n_k>1$, and exactly one of these has even length $k+1$.

To represent negative real numbers we apply the involution $x\mapsto S(x) = -1/x$.
This corresponds to the partition $\R\P^1=[0,\infty) \sqcup\, S\cdot [0,\infty)$.
Thus every $x\in\R\P^1\setminus \Q\P^1$ admits exactly one representation of the form $x=\Ecf{n_0,\dots}$ or $x=-1/\Ecf{n_0,\dots}$, and every $x\in \Q\P^1$ admits exactly two such representations of the same form, including $0=\Ecf{}=\Ecf{0}$ and $\infty = -1/\Ecf{} = -1/\Ecf{0}$.

Now consider the action of the modular group $\PSL_2(\Z)$ on $\R\P^1$ and of its euclidean submonoid generated by $L\& R$ on $[0,\infty]$.
If $x_i$ denotes the $i^{\mathrm{th}}$ remainder of $x\in [0,\infty]$, given by the tail $\Ecf{n_i,\dots}$ of its continued fraction expansion, then $x_0=(R^{n_0}L^{n_1})x_2$.
So the orbits of $x,y\in [0,\infty]$ under the euclidean monoid have non-empty intersection if and only if there exist even starting points $i,j$ at which the tails $x_i$ and $y_j$ coincide.
%
We deduce that $x,y\in \R\P^1$ belong to the same $\PSL_2(\Z)$-orbit if and only if there exist even starting points $i,j$ at which the tails $x_i$ and $y_j$ coincide.

\begin{Proposition}
\label{prop:Gauss-Galois-contfrac}
A geodesic in $\M$ has a lift $(\alpha_-,\alpha_+)\subset \HP$ whose endpoints satisfy $-1<\alpha_-<0$ and $1<\alpha_+<\infty$.
It intersects $\triangle$ along a sequence of triangles whose encoding in $\{L,R\}^{\Z}$ is obtained from the continued fraction expansions of $\alpha_+$ and $-1/\alpha_-$ by concatenating the transpose of the latter with the former.


The closed geodesics in $\M$ correspond to the periodic sequences, hence to the $\Gamma$-orbits of pairs $(\alpha_-,\alpha_+)$ of Galois conjugate quadratic irrationals.
\end{Proposition}

\subsection{Moduli space of elliptic curves and modular J-invariant}

A complex elliptic curve is isomorphic to the quotient of $\C$ by a lattice $\Lambda_\tau = \Z\oplus \tau \Z$ for a unique $\tau \in \HP$ up to the action of $\PSL_2(\Z)$.
Hence the modular orbifold corresponds to the moduli space of elliptic curves.

The field of meromorphic functions on $\C/\Lambda_\tau$, rather of $\Lambda_\tau$-invariant meromorphic functions on $\C$, is generated by the following Weierstrass sums over $\lambda \in \Lambda_\tau\setminus\{0\}$:
\begin{equation*}\textstyle
    \wp_\tau(u) = \sum' \frac{1}{(u-\lambda)^2}-\frac{1}{\lambda^2}
    \quad \mathrm{and} \quad
    \wp_\tau'(u) = -2\sum' \frac{1}{(u-\lambda)^3}
\end{equation*}
The Fourier expansion of $\wp_\tau(u)= u^{-2}+\sum_{n\ge1} (2n+1) G_{2n+2}(\tau) u^{2n}$ has coefficients the Eisenstein series $G_{n}(\Lambda)=\sum' \lambda^{-n}$ where the sum ranges over all $\lambda \in \Lambda_\tau \setminus \{0\}$, which are modular forms of weight $n$.

The functions $(x,y)=(\wp_\tau,\wp'_\tau)$ lie on the plane affine cubic $\E_\tau\colon y^2=4x^3-g_2x-g_3$ where $g_2(\tau)=60G_4(\tau)$ and $g_3(\tau)=140G_6(\tau)$, with discriminant $\Delta= g_2^3-27 g_3^2$.

The complex differential $1$-form $dx/y$ restricts on the cubic $\E_{\tau}$ to the holomorphic $1$-form $d\wp/\wp'=du$, whose integration along paths based at the origin $\infty$ lifts on its universal cover to the Jacobian integration map $\int^\infty du \colon H_1(\E_\tau;\R) \to \C$ sending its period lattice $H_1(\E_\tau;\Z)$ to its lattice of periods $\Z \omega_0+\Z\omega_1$ with $\omega_1/\omega_0=\tau$.

Hence the isomorphism class of $\C/\Lambda_\tau$ is equal to that of $\E_\tau$, characterised by the modular invariant $\J(\tau)=g_2^3(\tau)/\Delta(\tau)$. This $\PSL_2(\Z)$-modular holomorphic and bijective function $\J\colon \HP \to \C$ yields an isomorphism of complex curves $\M\simeq \C$.

Let us recall from \cite{Cox_primes-of-the-form_1997} that $\J(\tau)$ is algebraic if and only if $\tau$ is complex quadratic (in which case $\J(\tau)$ generates the maximal unramified abelian extension of the quadratic field $\Q(\tau)$ over that of $\Q$).
For example $\J(i)=12^3$, $\J(j)=0$.

\section{The modular torus and its mapping class group}

\subsection{The derived modular group and geodesics in the torus}

The abelianisation $\Z/2*\Z/3\to \Z/2\times \Z/3$ of the modular group corresponds to the $\Z/6$-Galois cover $\M' \to \M$ of the modular orbifold by the modular torus, which is a punctured torus whose fundamental group $\pi_1(\M')=\PSL_2(\Z)'$ is freely generated by 
\begin{equation*}
    X = [T^{-1},S] = LR =
    \begin{pmatrix}
    1 & 1\\
    1 & 2
    \end{pmatrix}
    \quad \mathrm{and} \quad
    Y = [T,S^{-1}] = RL =
    \begin{pmatrix}
    2 & 1\\
    1 & 1
    \end{pmatrix}
\end{equation*}
Note that $\Gamma'$ consists of all infinite order $A\in \Gamma$ with $\Rad(A)=0\bmod{6}$.

\begin{figure}[h]
    \centering
    \scalebox{0.5}{\begin{tikzpicture}[line cap=round,line join=round,>=triangle 45,x=4.166666666666667cm,y=4.166666666666667cm]
\clip(-1.2,-1.2) rectangle (1.2,1.2);

\fill[line width=0.pt,color=marron,fill=green,fill opacity=0.2] (-0.5,0.) -- (0.5,0.) -- (0.,1.) -- cycle;
\fill[line width=0.pt,color=marron,fill=red,fill opacity=0.2] (-0.5,0.) -- (0.5,0.) -- (0.,-1.) -- cycle;
\fill[line width=0.pt,color=marron,fill=red,fill opacity=0.2] (1.,0.) -- (0.5,0.) -- (0.,1.) -- cycle;
\fill[line width=0.pt,color=marron,fill=red,fill opacity=0.2] (-1.,0.) -- (-0.5,0.) -- (0.,1.) -- cycle;
\fill[line width=0.pt,color=marron,fill=green,fill opacity=0.2] (1.,0.) -- (0.5,0.) -- (0.,-1.) -- cycle;
\fill[line width=0.pt,color=marron,fill=green,fill opacity=0.2] (-1.,0.) -- (-0.5,0.) -- (0.,-1.) -- cycle;

\draw [line width=2.pt] (0.,0.) circle (4.166666666666667cm);

\draw [line width=2.pt,color=green] (1.,0.)-- (0.,1.);
\draw [line width=2.pt,color=green] (0.,1.)-- (-0.5,0.);
\draw [line width=2.pt,color=green] (-0.5,0.)-- (0.,-1.);
\draw [line width=2.pt,color=green] (0.,1.)-- (0.5,0.);
\draw [line width=2.pt,color=green] (0.5,0.)-- (0.,-1.);
\draw [line width=2.pt,color=green] (0.,1.)-- (-1.,0.);
\draw [line width=2.pt,color=green] (1.,0.)-- (0.,-1.);
\draw [line width=2.pt,color=green] (-0.5,0.)-- (-1.,0.);
\draw [line width=2.pt,color=green] (0.5,0.)-- (1.,0.);

\draw [line width=2.pt,color=marron] (0.5,0.)-- (0.666,0.333);
\draw [line width=2.pt,color=marron] (0.5,0.)-- (0.666,-0.333);
\draw [line width=2.pt,color=marron] (-0.5,0.)-- (-0.666,-0.333);
\draw [line width=2.pt,color=marron] (-0.5,0.)-- (-0.666,0.333);
\draw [line width=2.pt,color=marron] (0.,0.)-- (0.5,0.);
\draw [line width=2.pt,color=marron] (-0.5,0.)-- (0.,0.);

\draw [line width=2.5pt,color=forestgreen] (0.,1.)-- (1.,0.);
\draw [line width=2.5pt,color=forestgreen] (0.,1.)-- (0.,-1.);
\draw [line width=2.5pt,color=forestgreen] (-1.,0.)-- (0.,1.);
\draw [line width=2.5pt,color=forestgreen] (-1.,0.)-- (0.,-1.);
\draw [line width=2.5pt,color=forestgreen] (0.,-1.)-- (1.,0.);

\draw [line width=1.pt,color=black,-{Stealth[length=3.mm,width=2.5mm]}] (45:-0.65)-- (45:0.65);
\draw [line width=1.pt,color=black,-{Stealth[length=3.mm,width=2.5mm]}] (-45:-0.65)-- (-45:0.65);
\draw[color=black,anchor=south west] (45:0.7) node {\large$LR$};
\draw[color=black,anchor=north west] (-45:0.7) node {\large$RL$};

\begin{scriptsize}
\draw [fill=marron] (0.666,0.333) ++(-4.pt,0 pt) -- ++(4.pt,4.pt)--++(4.pt,-4.pt)--++(-4.pt,-4.pt)--++(-4.pt,4.pt);
\draw [fill=marron] (0.666,-0.333) ++(-4.pt,0 pt) -- ++(4.pt,4.pt)--++(4.pt,-4.pt)--++(-4.pt,-4.pt)--++(-4.pt,4.pt);
\draw [fill=marron] (-0.666,-0.333) ++(-4.pt,0 pt) -- ++(4.pt,4.pt)--++(4.pt,-4.pt)--++(-4.pt,-4.pt)--++(-4.pt,4.pt);
\draw [fill=marron] (-0.666,0.333) ++(-4.pt,0 pt) -- ++(4.pt,4.pt)--++(4.pt,-4.pt)--++(-4.pt,-4.pt)--++(-4.pt,4.pt);
\draw [fill=marron,rotate=45] (0.,0.) ++(-4.pt,0 pt) -- ++(4.pt,4.pt)--++(4.pt,-4.pt)--++(-4.pt,-4.pt)--++(-4.pt,4.pt);

\draw [fill=marron,shift={(0.5,0.)},rotate=270] (0,0) ++(0 pt,4.5pt) -- ++(3.8971143170299736pt,-6.75pt)--++(-7.794228634059947pt,0 pt) -- ++(3.8971143170299736pt,6.75pt);
\draw [fill=marron,shift={(-0.5,0.)},rotate=90] (0,0) ++(0 pt,4.5pt) -- ++(3.8971143170299736pt,-6.75pt)--++(-7.794228634059947pt,0 pt) -- ++(3.8971143170299736pt,6.75pt);

\draw [fill=black] (0.,1.) circle (2.5pt);
\draw[color=black] (0.,1.12) node {\large$\frac{1}{0}$};
\draw [fill=black] (1.,0.) circle (2.5pt);
\draw[color=black] (1.07,0.) node {\large$\frac{1}{1}$};
\draw [fill=black] (0.,-1.) circle (2.5pt);
\draw[color=black] (0.,-1.12) node {\large$\frac{0}{1}$};
\draw [fill=black] (-1.,0.) circle (2.5pt);
\draw[color=black] (-1.11,0.) node {\large$-\frac{1}{1}$};

\end{scriptsize}
\end{tikzpicture}}
    \scalebox{0.5}{\begin{tikzpicture}[line cap=round,line join=round,>=triangle 45,x=4.166666666666667cm,y=4.166666666666667cm]
\clip(-1.2,-1.2) rectangle (1.6,1.2);
\draw [line width=2.pt,color=black] plot[samples at={-150,-90,...,210},variable=\x] 
  (\x:1);

\fill[line width=0.pt,color=marron,fill=red,fill opacity=0.2] (0.,0.) -- (-150:1.) -- (-90:1.) -- cycle;
\fill[line width=0.pt,color=marron,fill=green,fill opacity=0.2] (0.,0.) -- (-90:1.) -- (-30:1.) -- cycle;
\fill[line width=0.pt,color=marron,fill=red,fill opacity=0.2] (0.,0.) -- (-30:1.) -- (30:1.) -- cycle;
\fill[line width=0.pt,color=marron,fill=red,fill opacity=0.2] (0.,0.) -- (150:1.) -- (90:1.) -- cycle;
\fill[line width=0.pt,color=marron,fill=green,fill opacity=0.2] (0.,0.) -- (90:1.) -- (30:1.) -- cycle;
\fill[line width=0.pt,color=marron,fill=green,fill opacity=0.2] (0.,0.) -- (150:1.) -- (-150:1.) -- cycle;

\draw [line width=2.pt,color=green] (0.,0.)-- (-150:1.);
\draw [line width=2.pt,color=green] (0.,0.)-- (-90:1.);
\draw [line width=2.pt,color=green] (0.,0.)-- (-30:1.);
\draw [line width=2.pt,color=green] (0.,0.)-- (30:1.);
\draw [line width=2.pt,color=green] (0.,0.)-- (90:1.);
\draw [line width=2.pt,color=green] (0.,0.)-- (150:1.);

\draw [line width=2.pt,color=forestgreen] (0.,0.)-- (-180:0.8657);
\draw [line width=2.pt,color=forestgreen] (0.,0.)-- (-120:0.8657);
\draw [line width=2.pt,color=forestgreen] (0.,0.)-- (-60:0.8657);
\draw [line width=2.pt,color=forestgreen] (0.,0.)-- (0:0.8657);
\draw [line width=2.pt,color=forestgreen] (0.,0.)-- (60:0.8657);
\draw [line width=2.pt,color=forestgreen] (0.,0.)-- (120:0.8657);

\draw [color=black,line width=1.pt, -{Stealth[length=3.mm,width=2.5mm]}] (0:1.1) to[out=90, in=-20,looseness=1] (60:1.);
\draw[color=black] (35:1.25) node {\large$ \Z/6$};

\draw[color=black] (0:1.35) node {\Huge$\rightsquigarrow$};

\begin{scriptsize}
\draw [line width=1.5pt,color=black, fill=white] (0:0) circle (3.5pt);
\draw [fill=marron, rotate around={45:(-180:0.8657)}] (-180:0.8657) ++(-4.pt,0 pt) -- ++(4.pt,4.pt)--++(4.pt,-4.pt)--++(-4.pt,-4.pt)--++(-4.pt,4.pt);
\draw [fill=marron, rotate around={10:(-120:0.8657)}] (-120:0.8657) ++(-4.pt,0 pt) -- ++(4.pt,4.pt)--++(4.pt,-4.pt)--++(-4.pt,-4.pt)--++(-4.pt,4.pt);
\draw [fill=marron, rotate around={-10:(-60:0.8657)}] (-60:0.8657) ++(-4.pt,0 pt) -- ++(4.pt,4.pt)--++(4.pt,-4.pt)--++(-4.pt,-4.pt)--++(-4.pt,4.pt);
\draw [fill=marron, rotate around={45:(0:0.8657)}] (0:0.8657) ++(-4.pt,0 pt) -- ++(4.pt,4.pt)--++(4.pt,-4.pt)--++(-4.pt,-4.pt)--++(-4.pt,4.pt);
\draw [fill=marron, rotate around={10:(60:0.8657)}] (60:0.8657) ++(-4.pt,0 pt) -- ++(4.pt,4.pt)--++(4.pt,-4.pt)--++(-4.pt,-4.pt)--++(-4.pt,4.pt);
\draw [fill=marron, rotate around={-10:(120:0.8657)}] (120:0.8657) ++(-4.pt,0 pt) -- ++(4.pt,4.pt)--++(4.pt,-4.pt)--++(-4.pt,-4.pt)--++(-4.pt,4.pt);

\draw [fill=marron,shift={(-150:1)},rotate=180] (0,0) ++(0 pt,4.5pt) -- ++(3.897pt,-6.75pt)--++(-7.794pt,0 pt) -- ++(3.897pt,6.75pt);
\draw [fill=marron,shift={(-90:1)},rotate=0] (0,0) ++(0 pt,4.5pt) -- ++(3.897pt,-6.75pt)--++(-7.794pt,0 pt) -- ++(3.897pt,6.75pt);
\draw [fill=marron,shift={(-30:1)},rotate=180] (0,0) ++(0 pt,4.5pt) -- ++(3.897pt,-6.75pt)--++(-7.794pt,0 pt) -- ++(3.897pt,6.75pt);
\draw [fill=marron,shift={(30:1)},rotate=0] (0,0) ++(0 pt,4.5pt) -- ++(3.897pt,-6.75pt)--++(-7.794pt,0 pt) -- ++(3.897pt,6.75pt);
\draw [fill=marron,shift={(90:1)},rotate=180] (0,0) ++(0 pt,4.5pt) -- ++(3.897pt,-6.75pt)--++(-7.794pt,0 pt) -- ++(3.897pt,6.75pt);
\draw [fill=marron,shift={(150:1)},rotate=0] (0,0) ++(0 pt,4.5pt) -- ++(3.897pt,-6.75pt)--++(-7.794pt,0 pt) -- ++(3.897pt,6.75pt);

\end{scriptsize}
\end{tikzpicture}
\begin{tikzpicture}[line cap=round,line join=round,>=triangle 45,
    rotate=270,
    x=6cm,y=6cm]
\clip(-0.2,-0.2) rectangle (1.3,0.6);
\draw [line width=2.pt,color=black] (0:0.8657) -- (30:1.);

\draw [line width=2.pt,color=green] (30.:0.05)-- (30:1.);

\draw [line width=2.pt,color=forestgreen] (0.:0.05)-- (0:0.8657);

\begin{scriptsize}
\draw [line width=1.5pt,color=black, fill=white] (0.05,0.015) ellipse (3pt and 4.5pt);
\draw [fill=marron, rotate around={45:(0:0.8657)}] (0:0.8657) ++(-4.pt,0 pt) -- ++(4.pt,4.pt)--++(4.pt,-4.pt)--++(-4.pt,-4.pt)--++(-4.pt,4.pt);

\draw [fill=marron,shift={(30:1)},rotate=0] (0,0) ++(0 pt,4.5pt) -- ++(3.897pt,-6.75pt)--++(-7.794pt,0 pt) -- ++(3.897pt,6.75pt);

\end{scriptsize}
\end{tikzpicture}}
    \caption{The free group $\PSL_2(\Z)'$ acts on $\HP$ with quotient a punctured torus $\M'$.
    The Galois group $\PSL_2(\Z)/\PSL_2(\Z)'=\Z/6$ acts on $\M'$ with quotient $\M$.}
    \label{fig:Hex-torus}
\end{figure}

Homotopy classes of loops in $\M'$ correspond to conjugacy classes in $\pi_1(\M')$, hence to reduced cyclic words in $X,X^{-1},Y,Y^{-1}$. 
The $n$-th power of the commutator $[X,Y]$ corresponds to the loop winding $n$ times around the puncture, and every other non-trivial homotopy class contains a unique geodesic. 

The cusp-compactification $\overline{\M'}=\Gamma'\backslash (\HP\sqcup \Q\P^1)$ of $\M'$ is homeomorphic to a torus, and has fundamental group $\pi_1(\overline{\M'})=\Gamma'/\triangleleft[X,Y]\triangleright=\Gamma'/\Gamma''=H_1(\Gamma';\Z)=\Z_X\oplus \Z_Y$. The simple loops and geodesics in $\M'$ correspond to those in $\overline{\M'}$.

The conjugacy classes in $\Gamma'$ associated to simple loops in $\M'$ correspond to the primitive vectors of the lattice $H_1(\Gamma';\Z)$.
The bases of $\Gamma'$ correspond to bases of $H_1(\Gamma';\Z)$, hence to pairs of simple loops in $\M'$ with one intersection point (their commutator in $\Gamma'$ corresponds to a loop circling once around the puncture).

The space of simple geodesics in $\M'$ corresponds to the spherisation $\S H_1(\M';\R)$.
%
The basis $(X,Y)$ of the lattice $H_1(\M';\Z)\subset H_1(\M';\R)$ identifies $\S H_1(\M';\Z)\subset \S H_1(\M';\R)$ with $\Q\S^1 \subset \R\S^1$.

\subsection{Mapping class group action on simple geodesics}

The mapping class group $\Mod(\overline{\M'})$ is identified to the outer automorphism group $\Out(\pi_1(\overline{\M'}))=\GL(\Z_X\oplus \Z_Y)$ by the \cite[Theorem 8.1]{Farb-Margalit_MappingClassGroups_2012} of Dehn-Nielsen-Baer,
and contains the orientation preserving mapping class group $\Mod^+(\overline{\M'}) = \SL_2(\Z_X \oplus \Z_Y)$ with index $2$ and quotient generated by $D_J\colon (X,Y) \mapsto (Y,X)$.

The group $\Mod^+(\overline{\M'}) = \SL(\Z_X\oplus \Z_Y)$ is generated by the positive Dehn twists $D_X$ and $D_Y$ along the simple loops $X$ and $Y$:
\begin{equation*}
    D_X\colon (X,Y)\mapsto (X,YX)
    \qquad 
    D_Y\colon (X,Y)\mapsto (XY,Y)
\end{equation*}
whose relations are generated by $D_YD_X^{-1}D_Y=D_X^{-1}D_YD_X^{-1}$ and $(D_YD_X^{-1})^6=1$.
It contains $D_X^{-1} D_Y D_X^{-1} = D_S \colon (X,Y)\mapsto (XYX^{-1},X^{-1})$, the positive rotation of order $4$ squaring to the elliptic involution $D_{-\Id} \colon(X,Y)\mapsto ((XY)X^{-1}(XY)^{-1},XY^{-1}X^{-1})$.

The group $\Mod(\M')$ of diffeotopy classes of diffeomorphisms of $\M'$ fixing the cusp fits into Birman's short exact sequence $1\to \pi_1(\M') \to \Mod(\M') \to \Mod(\overline{\M'})\to 1$, and the associated map $\Mod(\overline{\M'})\to \Out(\pi_1(\M'))$ is an isomorphism.
This restricts to $1\to \Z_X*\Z_Y \to \Map^+(\M') \to \SL_2(\Z_X\oplus \Z_Y)\to 1$. 
%

The action of $\Mod(\overline{\M'})=\GL_2(\Z)$ on $\S H_1(\M';\R) = \R\S^1$ quotients to an action of $\P\Mod(\overline{\M'})=\PGL_2(\Z)$ on $\P H_1(\M';\R) = \R\P^1$. 
(Beware not to confuse it with the action of $\Isom(\HP)=\PGL_2(\Z)$ on $\partial \HP = \R\P^1$: we will relate them later.)

We now use this action to determine for a primitive vector $(a,c)\in \Z^2$ the cyclic word in $X\& Y$ associated to the unique simple geodesic of $\M'$ homological to $X^aY^c$. 
After acting by $D_{S}$ we may assume $a,c\in \N^2$.
Expand $\frac{a}{c}=\Ecf{n_0,\dots, n_{k}}$ into a continued fraction of even length, and let $A = R^{n_0}\dots L^{n_{k}} \in \SL_2(\N)$ so that $A = \begin{psmallmatrix} a & b \\ c & d\end{psmallmatrix}$ sends $\frac{1}{0},\frac{0}{1}$ to its last two convergents $\frac{a}{c}, \frac{b}{d}$.
Then let $D_{A} = D_Y^{n_0}\dots D_X^{n_{k}}$ act by substitution on the seed $Y$ to yield the desired $X\& Y$-word, which is primitive.

This construction adapts to a geodesic of $\M'$ with any slope $(x,y) \in \S H^1(\M';\R)$ to yield its (dual) cutting sequence with the edges of the lattice $H^1(\M';\Z)=\Z_X\oplus \Z_Y$.
First act by $S$ to have $x,y >0$ and expand $\frac{x}{y}=\Ecf{n_0,n_1\dots}$ into a continued fraction. Then let $D_{(x,y)} = D_Y^{n_0}D_X^{n_{1}}\dots \in \SL_2(\N)$ act on the bi-infinite word $Y^\infty \in \{X,Y\}^{\Z}$.
For rational $(x,y) \in \S H_1(\M';\Z)$, this recovers the periodisation of the finite primitive word, or the dual cutting sequence of an affine translations of the associated line.
After replacing $X=LR$ and $Y=RL$, we recover the intersection pattern of the lifted geodesic $(\alpha_-,\alpha_+)\subset \HP$ with $\triangle$, whether the slope is rational or not. 

For $0/1 = \Ecf{}$, apply $D_{Y}^{0}$ to $Y^{\infty}$ to find $Y^\infty$, whose $L\& R$ conversion $(RL)^{\infty}$ corresponds to $(\tfrac{1}{2}(1-\sqrt{5}),\tfrac{1}{2}(1+\sqrt{5}))$.
For $1/1 = \Ecf{0,1}$, apply $D_{Y}^{0}D_{X}^{1}$ to $Y^{\infty}$ to find $(XY)^{\infty}$, whose $L\& R$ conversion $(LRRL)^{\infty}$ corresponds to $(-1-\sqrt{2},-1+\sqrt{2})$.
For $\operatorname{cotan}(2\pi/6)=\Ecf{0,1,\overline{1,2}}$, apply $(D_{Y}D_{X})(D_{Y}D_{X}^2)^\infty$ to $Y^\infty$... 

\begin{Proposition}[\cite{Series_Geo-Markov-Num_1985}]
\label{prop:Markov-Series-contfrac}
    Fix positive irrationals $-1/\alpha_-, \alpha_+$ and consider the $L\& R$-sequence obtained by concatenating their (transposed) continued fractions.
    
    The geodesic $(\alpha_-,\alpha_+)\subset \HP$ projects to a geodesic in $\M'$ which is simple if and only if the $L\& R$ can be paired to form an $X\& Y$-sequence which is Sturmian.
    
    The simple closed geodesics in $\M'$ correspond combinatorially to the periodic Sturmian $X\& Y$-sequences, and arithmetically to the Markov quadratic irrationals.
\end{Proposition}

\subsection{The Jacobian map and the Dedekind form}

The punctured hyperbolic torus $\M'$ is conformal to the complex elliptic curve $\overline{M'}$ punctured at its origin. The $\Z/6$-Galois symmetry of the cover $\overline{\M'}\to \overline{\M}$ implies that $\overline{M'}$ is isomorphic to the quotient $\C/\Z[j]$ of the complex plane by its hexagonal lattice underlying the ring of Eisenstein integers, which has invariant $\J(j)=0$.
%

The Jacobian map of $\overline{\M'}$ based at the projected cusp $\infty\in \partial \HP \mapsto \infty \in \overline{\M'}$ lifts to an identification $\int_\infty du \colon H_1(\M';\Z) \subset H_1(\M';\R) \to \Z\omega_0\oplus\Z\omega_2 \subset \C$ where the periods satisfy $\omega_2/\omega_0=\exp(i2\pi/3)$ and $\frac{1}{2}(\omega_0\wedge \omega_2)=\pi$ thus $\lvert \omega_0\rvert^2=4\pi/\sqrt{3}$.


Note that $\overline{\M'}$ is isomorphic over $\C$ to the affine cubic $y^2=4(x-j^0)(x-j^1)(x-j^2)$ with abelian differential $dx/y=-dx/\sqrt{4x^3-4}$ whose fundamental periods satisfy:
    \begin{equation*}
        \textstyle
        \frac{\varpi_k}{\varpi_0} = j^k
        \quad \mathrm{and}\quad 
        \varpi_0=2\int_\infty^1 \frac{-dx}{\sqrt{4x^3-4}} = \frac{\Upgamma(1/3)^3}{2^{4/3}\pi}
        \quad \mathrm{where} \quad
        \Upgamma(x)=\int_0^\infty t^x\exp(-t)\frac{dt}{t}.
    \end{equation*}
To be more precise, the complex curve $\M'$ with holomorphic form $du$ is isomorphic over $\R$ to the affine cubic $y^2=(x+i\lvert \omega_0\rvert j^0)(x+i\lvert \omega_0\rvert j^1)(x+i\lvert \omega_0\rvert j^2)$ with $dx/y$.

The lift on the universal cover $\HP\to \M'$ of the form $du$ is equal to $C\eta^4(\tau)d\tau$ where $\eta(\tau)=q^{1/24}\prod_{1}^\infty (1-q^n)$ is the Dedekind function, related to discriminant modular form by $\Delta=(2\pi)^{12}\eta^{24}$, and $C\in \C^*$.
Indeed, the cover $\M'\to \M$ has Galois group $\Z/6$ so $du^6$ descends on $\M$ and lifts on $\HP$ to a $\Gamma$-modular form of weight $12$, which must be a multiple of $\Delta(\tau)d\tau^6$.
We will compute $C$ later, but for now we may relate its modulus the half-area of $\M'$ by $\pi=\int\int_{\triangle(0,1,\infty)} \lvert C\eta^4(\tau)\rvert^2 \frac{d\tau\wedge d\overline{\tau}}{\Im(\tau)^2}$.

From the theta-series representation \cite[p. 30]{BGHZ_123-modular-forms_2008} of the function:
    \begin{equation*}
        \eta(\tau)=\sum_{n=0}^\infty \chi(n)q^{\tfrac{1}{24}n^2}
        \quad \text{we deduce} \quad
        \eta^4(\tau)=\sum_{n=0}^\infty \psi(n)q^{\tfrac{1}{24}n}
    \end{equation*}
where $\chi$ is the unique primitive character $\bmod{12}$ (that is the Jacobi symbol $\bmod{12}$, equal to $\chi(12n\pm 1)=1$, $\chi(12n\pm 5)=-1$ and $\chi(n)=0$ if $2\mid n$ or $3\mid n$), and where $\psi(n)$ is the sum of $\chi(abcd)$ over $\{(a,b,c,d)\in \N^4\mid a^2+b^2+c^2+d^2=n\}$.

The appendix \ref{sec:more-on-psi} records arithmetic properties of the function $\psi(n)$.
They rely on the fact that $\eta^4(6\tau)=q\prod_1^\infty (1-q^{6n})^4$ is \cite[\href{https://www.lmfdb.org/ModularForm/GL2/Q/holomorphic/36/2/a/a}{Modular Form 32.2.a.a}]{lmfdb}, the unique normalised cusp
eigenform for the group $\Gamma_0(36)= \{\begin{psmallmatrix}
a&b\\c&d \end{psmallmatrix}\in \Gamma \mid c\equiv 0 \bmod{36}\}$.

\section{The hexpunctured plane and its compactifications}

\subsection{Second derived modular group and hexpunctured plane}

The abelianisation $\Gamma'\to \Gamma'/\Gamma''$ associated to Hurwicz' map $\pi_1(\M') \to H_1(\M';\Z)$, given by $\Z_X*\Z_Y \to \Z_X\times \Z_Y$, corresponds to a Galois cover of $\M'$ by a lattice-punctured-plane $\M''$. 
The Jacobian integration map of $\M'$ based at the projected cusp $\infty\in \partial \HP \mapsto \infty \in \partial \M'$ lifts to an identification $\M'' \to H_1(\M';\R)\setminus H_1(\M';\Z)$.

The kernel $\Gamma''=\pi_1(\M'')$ is freely generated by an infinite set of elements (whose conjugacy classes in $\Gamma'$ are) indexed by $H_1(\M';\Z)=\Gamma'/\Gamma''$.
For example, the conjugates of $[X,Y]$ by all $X^mY^n\in \Gamma'$ with $m,n\in \Z$ freely generate $\Gamma''$.

We call $\hexp\colon \HP \to \M''$ the quotient by $\Gamma''$.
It maps the trivalent tree $\TT \subset \HP$ to a trivalent graph $\Hex \subset \M''$, and the vertices $\Gamma/\langle R\rangle\subset \partial \HP$ of $\triangle$ dual to the connected components of $\HP \setminus \TT$ to the lattice of cusps $\Gamma''\backslash \Gamma/\langle R\rangle$ in $\M''$.
Hence the inclusion $\Hex\subset \M''$ is a homotopy equivalence, and the regions in the complement $\M''\setminus \Hex$ are punctured regular hexagons (for the hyperbolic metric).

\begin{figure}[h]
    \centering
    \scalebox{0.9}{\begin{tikzpicture}[line cap=round,line join=round,>=triangle 45,x=2cm,y=2cm, rotate=-90]
\clip(-0.3,-1.2) rectangle (2.,1.2);
\draw [line width=1.5pt,color=black,shift={(0,0)}] plot[samples at={-150,-90,...,210},variable=\x] (\x:1);
\draw [line width=1.5pt,color=black,shift={(0+2*0.8657,0)}] plot[samples at={-150,-90,...,210},variable=\x] (\x:1);
\draw [line width=1.5pt,color=black,shift={(0.8657,2*0.75)}] plot[samples at={-150,-90,...,210},variable=\x] (\x:1);
\draw [line width=1.5pt,color=black,shift={(0.8657,-2*0.75)}] plot[samples at={-150,-90,...,210},variable=\x] (\x:1);

\draw [fill=marron, rotate around={45:(0:0.8657)}] (0:0.8657) ++(-4.pt,0 pt) -- ++(4.pt,4.pt)--++(4.pt,-4.pt)--++(-4.pt,-4.pt)--++(-4.pt,4.pt);

\draw [fill=marron,shift={(30:1)},rotate=0] (0,0) ++(0 pt,4.5pt) -- ++(3.897pt,-6.75pt)--++(-7.794pt,0 pt) -- ++(3.897pt,6.75pt);

\draw [color=grey,line width=0.5pt, -{Stealth[length=1.5mm,width=1.mm]},shift={(30:1)}] (-80:0.2) to[out=0, in=-60,looseness=1] (20:0.2);
\draw[color=grey,shift={(30:1)}, anchor=west] (-35:0.45) node {$T^{-1}$};

\draw [color=grey,line width=0.5pt, -{Stealth[length=1.5mm,width=1.mm]},shift={(30:1)}] (-100:0.2) to[out=180, in=-120,looseness=1] (160:0.2);
\draw[color=grey,shift={(30:1)}, anchor=west] (-145:0.45) node {$T$};

\draw [color=grey,line width=0.5pt, -{Stealth[length=1.5mm,width=1.mm]},shift={(0:0.8657)}] (-80:0.15) to[out=0, in=0,looseness=1.5] (80:0.15);
\draw[color=grey,shift={(0:0.8657)}] (0:0.39) node {$S^{-1}$};

\draw [color=grey,line width=0.5pt, -{Stealth[length=1.5mm,width=1.mm]},shift={(0:0.8657)}] (-100:0.15) to[out=180, in=180,looseness=1.5] (100:0.15);
\draw[color=grey,shift={(0:0.8657)}] (-175:0.3) node {$S$};

\draw [color=red,line width=0.5pt, -{Stealth[length=1.5mm,width=1.mm]}] (-40:0.5) to[out=80, in=-20,looseness=1] (80:0.4);
\draw[color=red] (0:0.1) node {$R$};

\draw [color=blue,line width=0.5pt, -{Stealth[length=1.5mm,width=1.mm]},shift={(2*0.8657+0.05,0)}] (-130:0.4) to[out=110, in=-140,looseness=1] (110:0.4);
\draw[color=blue,shift={(2*0.8657-0.1,0)}] (0:0.1) node {$L$};
\end{tikzpicture}}
    \hspace{1cm}
    \scalebox{0.9}{\begin{tikzpicture}[rotate=30]
\clip[shift={(2.5,2.5)}] (-30:1) circle (26mm);

\begin{scope}[shift={(2*2*0.8657,2*2*0.75)}]
    \draw[color=black, line width=3pt] (-120:0.8657) -- (-90:1);
    \draw[color=black, line width=3pt, anchor=north] (-120:0.8657) node {$i$};
    \draw[color=black, line width=3pt, anchor=west] (-90:1) node {$j$};
\end{scope}

\foreach \i in {0,...,3}
\foreach \j in {1,3} {
\begin{scope}[shift={(2*\i*0.8657+0.8657,2*\j*0.75)}]
\draw [line width=1.5pt,color=black] plot[samples at={-150,-90,...,210},variable=\x] (\x:1);
\draw [fill=marron, rotate around={45:(0:0.8657)}] (0:0.8657) ++(-3.pt,0 pt) -- ++(3.pt,3.pt)--++(3.pt,-3.pt)--++(-3.pt,-3.pt)--++(-3.pt,3.pt);
\draw [fill=marron, rotate around={45:(-180:0.8657)}] (-180:0.8657) ++(-3.pt,0 pt) -- ++(3.pt,3.pt)--++(3.pt,-3.pt)--++(-3.pt,-3.pt)--++(-3.pt,3.pt);
\draw [line width=0.8pt,color=black] (0,0) circle (2.5pt);
\end{scope}
}

\foreach \i in {0,...,3}
\foreach \j in {0,2} {
    \begin{scope}[shift={(2*\i*0.8657,2*\j*0.75)}]
        \draw [line width=1.5pt,color=black] plot[samples at={-150,-90,...,210},variable=\x] (\x:1);
        \draw [fill=marron, rotate around={45:(-180:0.8657)}] (-180:0.8657) ++(-3.pt,0 pt) -- ++(3.pt,3.pt)--++(3.pt,-3.pt)--++(-3.pt,-3.pt)--++(-3.pt,3.pt);
        \draw [fill=marron, rotate around={15:(-120:0.8657)}] (-120:0.8657) ++(-3.pt,0 pt) -- ++(3.pt,3.pt)--++(3.pt,-3.pt)--++(-3.pt,-3.pt)--++(-3.pt,3.pt);
        \draw [fill=marron, rotate around={-15:(-60:0.8657)}] (-60:0.8657) ++(-3.pt,0 pt) -- ++(3.pt,3.pt)--++(3.pt,-3.pt)--++(-3.pt,-3.pt)--++(-3.pt,3.pt);
        \draw [fill=marron, rotate around={45:(0:0.8657)}] (0:0.8657) ++(-3.pt,0 pt) -- ++(3.pt,3.pt)--++(3.pt,-3.pt)--++(-3.pt,-3.pt)--++(-3.pt,3.pt);
        \draw [fill=marron, rotate around={15:(60:0.8657)}] (60:0.8657) ++(-3.pt,0 pt) -- ++(3.pt,3.pt)--++(3.pt,-3.pt)--++(-3.pt,-3.pt)--++(-3.pt,3.pt);
        \draw [fill=marron, rotate around={-15:(120:0.8657)}] (120:0.8657) ++(-3.pt,0 pt) -- ++(3.pt,3.pt)--++(3.pt,-3.pt)--++(-3.pt,-3.pt)--++(-3.pt,3.pt);
        \draw [fill=marron,shift={(-150:1)},rotate=180] (0,0) ++(0 pt,{4.5/1.5pt}) -- ++({3.897/1.5pt},{-6.75/1.5pt})--++({-7.794/1.5pt},0 pt) -- ++({3.897/1.5pt},{6.75/1.5pt});
        \draw [fill=marron,shift={(-90:1)},rotate=0] (0,0) ++(0 pt,{4.5/1.5pt}) -- ++({3.897/1.5pt},{-6.75/1.5pt})--++({-7.794/1.5pt},0 pt) -- ++({3.897/1.5pt},{6.75/1.5pt});
        \draw [fill=marron,shift={(-30:1)},rotate=180] (0,0) ++(0 pt,{4.5/1.5pt}) -- ++({3.897/1.5pt},{-6.75/1.5pt})--++({-7.794/1.5pt},0 pt) -- ++({3.897/1.5pt},{6.75/1.5pt});
        \draw [fill=marron,shift={(30:1)},rotate=0] (0,0) ++(0 pt,{4.5/1.5pt}) -- ++({3.897/1.5pt},{-6.75/1.5pt})--++({-7.794/1.5pt},0 pt) -- ++({3.897/1.5pt},{6.75/1.5pt});
        \draw [fill=marron,shift={(90:1)},rotate=180] (0,0) ++(0 pt,{4.5/1.5pt}) -- ++({3.897/1.5pt},{-6.75/1.5pt})--++({-7.794/1.5pt},0 pt) -- ++({3.897/1.5pt},{6.75/1.5pt});
        \draw [fill=marron,shift={(150:1)},rotate=0] (0,0) ++(0 pt,{4.5/1.5pt}) -- ++({3.897/1.5pt},{-6.75/1.5pt})--++({-7.794/1.5pt},0 pt) -- ++({3.897/1.5pt},{6.75/1.5pt});
        \draw [line width=0.8pt,color=black] (0,0) circle (2.5pt);
    \end{scope}
    }
    

\draw[color=black,anchor=east] (2*1*0.8657,2*2*0.75) node{$-1$};
\draw[color=black,anchor=west] (2*2.5*0.8657,2*1*0.75) node{$+1$};
\draw[color=black,anchor=south] (2*2*0.8657,2*2*0.75) node {$\infty$};
\draw[color=black,anchor=north] (2*1*0.8657+0.8657,2*1*0.75) node {$0$};

\begin{scope}[shift={(2*2*0.8657,2*2*0.75)}]
    \draw [color=black,line width=1.5pt, -{Stealth[length=2.mm,width=1.5mm]}, anchor=south] (-120:0.8657) -- (-30:2*0.75) node {$RL$};
    \draw [color=black,line width=1.5pt, -{Stealth[length=2.mm,width=1.5mm]}, anchor = north] (-120:0.8657) -- (-79:2.3) node{$LR$};
\end{scope}
\end{tikzpicture}}
    \caption*{Action of $\Gamma/\Gamma''$ on the hexagonal graph $\Hex$ and cusps $\Gamma''\backslash \Gamma /\langle R\rangle$.}
\end{figure}

We now describe the Galois action of $\Gamma/\Gamma''$ on the metabelian cover $\M'' \to \M$ combining the action of $\Gamma'/\Gamma''$ on $\M'' \to \M'$ with the action of $\Gamma/\Gamma'$ on $\M'\to \M$.

The group $\Gamma/\Gamma''$ acts freely transitively on the oriented edges of $\Hex$. 
More precisely the generators $S$ and $T$ act by rotations of order $2$ and $3$ around the mid-edges and vertices of $\Hex$.
The subgroup $\Gamma'/\Gamma''$ generated by $X=LR$ and $Y=RL$ acts by translation of $\Hex$ with fundamental domain a tripod formed by three adjacent edges. 

Consequently the action of $A\in \Gamma$ on $\Hex$ sends the base edge to another edge whose angle is determined by $A\bmod{\Gamma'}$ and whose base point is determined by $A^6\bmod{\Gamma''}$. 
In particular $\Gamma'/\Gamma''$ identifies with the edges of $\Hex$ that are parallel to the base edge.

\begin{Proposition}[\cite{CLS_phdthesis_2022}]
\label{prop:CLS_hexagonal-group}
The group $\Gamma/\Gamma''$ is the semi-direct product $\Gamma'/\Gamma'' \rtimes \Gamma/\Gamma'$
where the action of the quotient $\Gamma/\Gamma'=\Z/2\times \Z/3$ by outer-automophisms on the kernel $\Gamma'/\Gamma''=\Z_X\oplus \Z_Y$ is given by $\AD_S=D_{S^2}=D_{-\Id}$ and $\AD_T=D_{T^2}$ since we have $SXS^{-1}=X^{-1}$, $SYS^{-1}=Y^{-1}$ and $TXT^{-1}=Y^{-1}$, $TYT^{-1}=-XY^{-1}$. 

This represents $\Gamma/\Gamma''$ as the affine isometry group of the oriented hexagonal lattice $H_1(\Gamma';\Z)=\Z_X\oplus \Z_Y$ with $angle{(X,Y)}=\frac{2\pi}{6}$, where the translation action of $\Gamma'/\Gamma''$ is by $\Z^2$-translation while the outer-automorphism action of $\Gamma/\Gamma'$ is by $\Z/6$-rotation.    
\end{Proposition}

\begin{Remark}[Uniformization]
Let $\Lambda = \omega_0 \Z[j]$ be the hexagonal period-lattice of $\eta^4$, that is $\partial \hexp(\Q\P^1) \subset \C$, and denote $\wp_\Lambda$ its associated Weierstrass function.

The functions $\wp_\Lambda \circ \hexp \colon \HP \to \C$ and $\wp'_\Lambda \circ \hexp \colon \HP \to \C$ are entire, $\Gamma'$-invariant and $\Gamma/\Gamma'$-equivariant.
Their poles on the boundary are precisely $\Q\P^1$.
Together, they define the map $\HP \to \C^2$ uniformizing the affine elliptic curve $y^2=4(x^3-\omega_0^3)$.  
\end{Remark}

\subsection{Cusp-compactification of the hexponential map}

The map $\hexp \colon \HP \to \M''$ in the coordinates $\{z\in \C \mid \Im(z)>0\} \to \C\setminus \omega_0 \Z[j]$ is the primitive of the function $C\eta(z)^4dz$, with expansion:
\begin{equation}\label{eq:hexp}\tag{hexp}
    \hexp(\tau)=\int_\infty^\tau C\eta^4(z)dz = \frac{12C}{i\pi}\sum_{n=1}^\infty \frac{\psi(n)}{n}\cdot q^{\tfrac{1}{24}n}
\end{equation}
The fundamental period of $\M'$ is $\omega_0 = \frac{1}{i}\lvert \omega_0\rvert = \int_\infty^0 C\eta^4(\tau)d\tau$, and by \cite{Shalev_eta-period-stack-exchange_2019} we have:
\begin{equation*}\textstyle
    \int_\infty^0 \eta^4(\tau)d\tau = \frac{1}{i} \frac{3^{1/2}}{2^{7/3}\pi^2} \Upgamma(1/3)^3
    \quad \mathrm{thus} \quad
    C = \frac{2^{10/3}}{3^{3/4}} \frac{\pi^{5/2}}{\Upgamma(1/3)^{3}}
\end{equation*}
We deduce from the area of $\M'$ the value of $\int\int_{\triangle(0,1,\infty)} \lvert \eta^4(\tau)\rvert^2 \frac{d\tau\wedge d\overline{\tau}}{\Im(\tau)^2}=\frac{\pi}{\lvert C\rvert^2}$.

Observe that the lattice relation $\partial \hexp(0)=2\hexp(i)$ expresses the value of an conditionally convergent series and of an absolutely convergent series:
\begin{equation*}\textstyle
    \sum_{1}^\infty \tfrac{\psi(n)}{n}=
    2\sum_{1}^\infty\frac{\psi(n)}{n} \exp\left(\frac{-\pi}{12}n\right) 
    \quad \text{as equal to} \quad
    \frac{\pi\lvert\omega_0\rvert}{12C}=
    \frac{1}{24}
    \frac{3^{1/2}}{2^{1/3}} \frac{\Upgamma(1/3)^{3}}{2\pi}.
\end{equation*}
Let us generalise the computation of such special values to obtain Theorem \ref{thm:values-hexp(rational)}, whose proof will justify the conditional convergence. 

Recall that the cusp map $\partial \hexp \colon \Gamma/\langle R \rangle \to \Gamma''\backslash \Gamma/\langle R \rangle$ sends the vertices $\Q\P^1\subset \triangle$ to the cusps $H_1(\M';\Z) \subset H_1(\M';\R)$.
The concatenation of even continued fraction expansions of rational numbers quotients $\bmod{\Gamma''}$ to the addition of cusps twisted by the $\Gamma/\Gamma'$-action (so the twist is trivial if the first rational number belongs to $\Gamma'\cdot \infty$).
This yields twisted-addition identities between the values of the $\partial \hexp$ at various $\alpha\in \Q\P^1$, given by the improper integral and conditionally convergent series \eqref{eq:hexp}.

In particular, for every $A=\begin{psmallmatrix} a & b \\ c & d \end{psmallmatrix}\in \Gamma'$ equal to $X^mY^n \mod{\Gamma''}$, we have 
\begin{align*}
    \partial\hexp(a/c)&= \lvert \omega_0\rvert\left(m\exp(-\tfrac{i\pi}{6})+n\exp(+\tfrac{i\pi}{6})\right) \\
    \partial\hexp(b/d)&= \lvert \omega_0\rvert\left(m\exp(-\tfrac{i\pi}{6})+n\exp(+\tfrac{i\pi}{6})-i\right)
\end{align*}
since the vectors $X,Y$ act on on $\M''=\C\setminus (\omega_0\Z[j])$ by translations of $\lvert \omega_0\rvert\exp(\pm \tfrac{i\pi}{6})$.
Note that the left-columns of matrices in $\Gamma'$ yield all rationals.

\begin{Theorem}[Values of \texorpdfstring{$\partial \hexp$}{dhexp}]
\label{thm:values-hexp(rational)}
    Consider $a/c\in \Q\P^1$ with $\gcd(a,b)=1$.
    
    There exists $A=\begin{psmallmatrix} a & b \\ c & d \end{psmallmatrix} \in \Gamma'$, and we may consider its class $A\equiv X^mY^n \mod{\Gamma''}$. 
    The value $\partial \hexp(\frac{a}{c})$ of the improper integral and conditionally convergent series:
    \begin{equation*}
        \partial \hexp(\tfrac{a}{c}) = 
        \int_\infty^{\frac{a}{c}} C\eta^4(z)dz 
        = \tfrac{12C}{i\pi}\sum_{n=1}^\infty \tfrac{\psi(n)}{n}\cdot \exp\left(\tfrac{i\pi}{12}\tfrac{a}{c}n\right)
        \quad \text{where} \quad
        C = \tfrac{2^{10/3}}{3^{3/4}} \tfrac{\pi^{5/2}}{ \Upgamma(1/3)^{3}} 
    \end{equation*}
    \begin{equation*}
        \text{is equal to}\quad 
        \partial \hexp(\tfrac{a}{c}) = 
        \lvert \omega_0 \rvert \left(m\exp(-\tfrac{i\pi}{6})+n\exp(+\tfrac{i\pi}{6})\right)
        \quad \text{where} \quad
        \lvert \omega_0 \rvert = \tfrac{2\pi^{1/2}}{3^{1/4}}.
    \end{equation*}
\end{Theorem}

\begin{proof}
We are left with the conditional convergence of the series $\sum \tfrac{\psi(n)}{n} e^{i2\pi r}$ for $r\in \Q$.

Let us first explain why they are not absolutely convergent. 

The values $n\in \N$ such that $\psi(n)\ne 0$ are characterised in \cite[Section 2.2]{Serre_lacunarite-puissances-eta_1985}, in terms of the ramification of $n$ in the ring of integers of the quadratic field $\Q(j)$.
In particular $\psi(p)\ne 0$ for primes $p\equiv 1\bmod{6}$. 
Moreover, the Dirichlet's density theorem \cite[§4.1]{Serre_arithmetique_1970} implies that the sum $\sum \tfrac{1}{p}$ over primes $p\equiv 1\bmod{6}$ diverges.


Now let us explain the convergence for $r=\tfrac{a}{24c}$ for $\gcd(a,c)=1$.

Recall thath $\psi(n)$ are the Fourier coefficients of $\eta^4(6\tau)$ which is a cusp form of weight $2$ for the group $\Gamma_0(36)$.
We may thus apply \cite[Corollary 1.5.21]{Jutila_Lectures-exponential-sums_1987}, saying that for all $\epsilon>0$ we have $\lvert \sum_1^N \psi(n)\exp(i2\pi r)\rvert = O(c^{2/3}N^{5/6+\epsilon})$. Hence the $L$-series $L(s,r)=\sum_1^\infty \frac{\psi(n)}{n^s}\exp(i2\pi r)$ converges where $\Re(s)\ge 5/6$, in particular at $s=1$.
\end{proof}

\begin{Question}[Convergence in $\Q_p$]
Does the series defined by $\partial \hexp(r)$ for $r\in \Q$ converge in other $p$-adic completions of $\Q$?

In this direction, note that $\hexp$ is a non-congruence modular function, namely for the group $\Gamma''$, so by \cite{Li-Long_Fourier-coefficient--noncongruence-cuspforms_2012} the denominators of its coefficients $\tfrac{1}{n}\psi(n)$ are unbounded.

Besides, the work of \cite{Pollack-Stevens_Overconvergent-modular-symbols_2011} on overconvergent modular symbols and $p$-adic $L$-functions may be applied to the group $\Gamma'$.
\end{Question}

\newpage

\subsection{Spheric compactification of the hexpunctured plane}
\label{subsec:Shexp}
%
A complete geodesic in $\HP$ or $\TT$ encoded by a bi-infinite $L\& R$-word maps to a complete geodesic in $\M''$ or $\Hex$ encoded by the same bi-infinite $L\& R$-word.
Observe that a euclidean straight line in $\M''$ determines a unique complete geodesic in $\M''$ both topologically (they are a homotopic relative to endpoints) and combinatorially (they intersect the triangulation $\hexp(\triangle)$ dual to $\Hex$ in the same bi-infinite sequence of $\Gamma$-translates of the triangle $\hexp(0,1,\infty)$).

The left translation action of $\Gamma$ on $\TT\subset \HP$ quotients to left translation-rotation action of $\Gamma/\Gamma''$ on $\Hex\subset \M''$.
Hence the primitive conjugacy classes of $\Gamma$ correspond to the bi-infinite periodic geodesics of $\Hex$ passing through the oriented base edge.
One may construct various geometricombinatorial invariants \cite[Section 3.2]{CLS_phdthesis_2022}.
For instance the Rademacher number of a hyperbolic element in $\Gamma$ can be expressed as the total asymptotic winding number of its axis projected in $\Hex\subset \M''$.

We define a map $\Shexp \colon \Radial \to \S H_1(\M';\R)$ from the subset of numbers $\alpha \in \R\P^1$ such that as $\tau \in \HP$ converges to $\alpha \in \partial \HP$, the modulus $\lvert \hexp(\tau) \rvert \in \R_+$ diverges while the argument $\arg \hexp(\tau) \in \R/(2\pi \Z)$ converges to a limit defining $\Shexp(\alpha)$.

Thus $\Radial \subset \R \setminus \Q$ contains the numbers $\alpha$ such that the geodesic $(i,\alpha) \subset \HP$ maps under $\hexp$ to a geodesic of $\M''$ which escapes in a definite direction (it follows the geodesic of $\Hex$ encoded by the $L\& R$-continued fraction expansion of $\alpha$).

In other terms, an irrational $\alpha = \Ecf{n_0, n_1, \dots}$ belongs to $\Radial$ if and only if the orbit of $\infty \in\Gamma''\backslash \Gamma/ \langle R \rangle = H_1(\M';\Z)$ under the $A_k = R^{n_0}\dots  L^{n_{2k+1}}\in \Gamma/\Gamma''$ projectifies to a Cauchy sequence in $\S(\Gamma''/\Gamma')=\S H_1(\M';\Z)$.

The surjective maps $\Shexp \colon \Radial \to \S H_1(\M';\R)$ and $\partial \hexp \colon \Q\P^1 \to H_1(\M';\Z)$ have disjoint domains. Together, they define (for obvious topologies) a continuous compactification of $\hexp \colon \HP \to \M''$, equivariant under the left actions of $\Gamma$ and $\Gamma/\Gamma''$.
The translation action of $\Gamma'$ on $\S H_1(\M';\R)$ is trivial, so we have an equivariant map $\Shexp \colon \Radial \bmod{\Gamma'} \to \S H_1(\M';\R)$ under the actions of $\Gamma$ and $\Gamma/\Gamma'$.


\begin{Remark}[Singularities along the boundary]

If $\alpha \in \R\P^1$ is neither rational nor quadratic then $\hexp(\tau)$ diverges as $\tau\in \HP$ diverges to $\alpha\in \R\P^1$.
The subset $\Radial \subset \R\P^1$ is $\Gamma$-invariant and dense in $\R\P^1$.
The map $\Shexp \colon \Radial \to \S H_1(\M',\R)$ is uniformly continuous and surjective.
Suppose $\alpha \in \R\P^1$ is a fixed point of $A\in \Gamma$: if $A^6\notin \Gamma''$ then $\alpha \in \Radial$, but if $A^6\in \Gamma''$ then $\alpha \notin \Radial$ yet $\lim_\alpha \hexp(\tau) = \partial\hexp(r)$ for some $r\in \Q\P^1$.
\end{Remark}

\begin{proof}
    The proof follows from the combinatorial and topological behaviour of $\hexp$, using the fact that the symbolic encoding of real numbers by continued fractions is surjective and uniformly continuous (if $\alpha,\beta\in \R$ have continued fractions with the same prefix given by $A=R^{n_0}\cdots L^{n_{2k+1}}$ for some $k\in \N$, then $\lvert \alpha-\beta\rvert \le 2/(F_{2k}F_{2k+1})$ where $(F_k)_{k\in \N}$ is the usual Fibonacci sequence). 
\end{proof}

\section{Transcendental properties of the hexponential}

\subsection{The Inshection map and Sturmian numbers}

We defined a surjective map $\Shexp \colon \Radial \bmod{\Gamma'} \to \S H_1(\M';\R)$ where $\Radial \subset \R\P^1\setminus \Q\P^1$, sending uncountably many distinct $\Gamma'$-orbits in $\Radial$ to the same slope in $\S H_1(\M';\R)$.
Let us define a section $\InSh \colon \S H_1(\M';\R) \to \Sturm \bmod{\Gamma'}$, which is $\Gamma/\Gamma'$-equivariant, injective and whose image will have ``index at most $2$'' in the set $\Sturm$ to be defined.


A slope $(x,y) \in \S H_1(\M';\R)$ corresponds to a complete simple geodesic $\alpha\subset \M'$: it lifts in $\HP$ to the $\Gamma'$-orbit of a geodesic $(\alpha_-,\alpha_+)$ and we let $\InSh(x,y)=\Gamma'\cdot \alpha_+$.
The geodesic $\alpha \subset \M'$ lifts in $\M''$ to the $\Gamma'/\Gamma''$-orbit of geodesics which intersect the triangulation $\hexp(\triangle)$ according to the $X^{\pm1}\& Y^{\pm1}$ dual cutting sequences of the euclidean lines inside $H_1(\M';\R)\setminus H_1(\M';\Z)$ of slope $(x,y)$.
The $X\& Y$ cutting sequences of such lines are precisely the Sturmian sequences, and two sequences corresponds to a same $\Gamma'/\Gamma''$-orbit if and only if they differ by a shift.
Replacing $X=LR$, $Y=RL$ recovers the bi-infinite $L\& R$-sequences giving the continued fraction expansions of the numbers $-1/\alpha_-,\alpha_+$.

Let us remark that in $\M''$, the lifted geodesics $(\alpha_-,\alpha_+)\bmod{\Gamma''}$ intersect the cusp-to-cusp geodesic $(\infty,0)\bmod{\Gamma''}$ a finite number of times if and only if $x/y$ is rational, and they approach the cusps arbitrarily close if and only if the terms in the continued fraction expansion of $x/y$ are unbounded.


We call $\Sturm$ the set of \emph{Sturmian} numbers $\alpha_+,-1/\alpha_-$ obtained as above, namely whose continued fraction expansions arise from Sturmian $X^{\pm1}\& Y^{\pm1}$-sequences by replacing $X=RL$ and $Y=RL$.
It is endowed with the sign reversing involutions $\sigma \colon \alpha_-\leftrightarrow \alpha_+$ and $S\colon \alpha \mapsto -1/\alpha$, which commute and compose to $\epsilon \colon \alpha_+ \mapsto -1/\alpha_-$.

Note that $\hexp$ is also invariant under $\epsilon$, namely $\hexp(\alpha_+)=\hexp(-1/\alpha_-)$, and the map $\S H_1(\M';\R) \to \Sturm \bmod{\{\Gamma', \epsilon\}}$ is bijective. In this sense the image of $\InSh$ has index at most $2$ in $\Sturm \bmod{\Gamma'}$, taking into account the fact that some $\alpha\in \Sturm$ satisfy $\alpha=\epsilon(\alpha) \bmod{\Gamma'}$.
Moreover $\InSh$ is $\Gamma/\Gamma'$-equivariant in the sense that 
\begin{align*}
    &\InSh \AD_S(x,y)=S(\InSh(\theta))
    \quad \mathrm{namely} \quad 
    \InSh(-x,-y)=-1/\InSh(\theta)
    \\
    &\InSh \AD_T(x,y)=T(\InSh(\theta))
    \quad \mathrm{namely} \quad
    \InSh(-y,x-y)=(\InSh(\theta)-1)/\InSh(\theta)
\end{align*} 
so $\InSh$ descends to a map $\S H_1(\M';\R) \bmod{\AD_\Gamma} \to \Sturm\bmod{\Gamma}$.


\begin{Remark}[Cantor set]
    The set $\Sturm \subset \R\P^1$ is compact, with a dense countably infinite set of isolated points whose complement is a Cantor set. By \cite{Pelczynski_zero-dimensional-spaces_1965} such properties describe a unique space up to homeomorphism.
    It has Hausdorff dimension zero.
    The map $\Shexp \colon \Sturm \to \S H_1(\M';\R)$ is continuous and surjective.
\end{Remark}

\subsection{Diophantine approximation of Sturmian numbers}


For a real $\alpha \in \R$, its Lagrange constant $\Lag(\alpha)$ is the supremum of $L\in \R_+$ such that there are infinitely many $p,q\in \Z\times \N^*$ with $\lvert \alpha -p/q \rvert < 1/(L q^2)$, namely:
\begin{equation*}\textstyle
    \tfrac{1}{\Lag(\alpha)}= \liminf_{q}({q^2\lvert \alpha -p/q \rvert})
\end{equation*}
It may be expressed in terms of the continued fraction expansion $\alpha = \Ecf{a_0,a_1,\dots}$:
\begin{equation*}\textstyle
    \Lag(\alpha)= \limsup_{n}({\Ecf{0, a_{n-1},\dots, a_{0}}+\Ecf{a_{n},a_{n+1},\dots}})
\end{equation*}
In particular $\Lag(\alpha)=0 \iff \alpha \in \Q$ whereas $\Lag(\alpha)=+\infty \iff \limsup_n(a_n) = +\infty$.
Moreover, for $\alpha \notin \Q$ we have $\Lag(\alpha)\ge \limsup_n \left(\Ecf{0, 1,\dots, 1}+\Ecf{1,1,\dots}\right) = \sqrt{5}$.

In $\M$, the horoball neighborhood $\B(h)$ of the cusp $\infty$ at height $h\ge 1$ is the image of the projection $\{\Im(z) > h\} \subset \HP$, which has area $1/h$.
We will also consider its lift in $\M'$ which has area $6/h$, as the cover $\M'\to \M$ restricts to a degree $6$ cover of punctured discs between these cusp neighborhoods.

The geodesic $(\infty,\alpha)\in \HP$ projects to a geodesic in $\M'$ whose penetrations into $\B(h)$ are indexed by the $n\in \N$ such that $\tfrac{1}{2}\left(\Ecf{0, a_{n-1},\dots, a_{0}}+\Ecf{a_{n},a_{n+1},\dots}\right)>h$.

For distinct $\alpha_-,\alpha_+ \in \R\P^1$, we define the Markov-Cohn constant $\Cohn(\alpha_-,\alpha_+)$ as the supremum of $2h\in \R_+$ such that $\B(h)$ is disjoint from the geodesic $(\alpha_-,\alpha_+)\subset \M$.

When $-1< \alpha_- \le 0$ and $1 \le \alpha_+ < \infty$ have continued fraction expansions $-1/\alpha_- = \Ecf{a_{-1}, a_{-2}, \dots}$ and $\alpha_+=\Ecf{a_0, a_{1}, a_{2}, \dots}$, the geodesic $(\alpha_-,\alpha_+)\subset \HP$ intersects $\triangle$ according to the sequence of $a_{n}$ for $n\in \Z$.
It projects to a geodesic in $\M$ which penetrates $\B(h)$ each time $\tfrac{1}{2}\left(\Ecf{0, a_{n-1}, a_{n-2}, \dots}+\Ecf{a_{n}, a_{n+1}, \dots}\right)>h$, thus
\begin{equation*}\textstyle
    \Cohn(\alpha_-,\alpha_+)=\sup_n (\Ecf{0, a_{n-1}, a_{n-2}, \dots}+\Ecf{a_{n}, a_{n+1}, \dots})
\end{equation*}
The following proposition is due to Haas \cite{Haas_geometry-Markoff-forms_1987} building on work of Cohn \cite{Cohn_Markoff-perforated-torus_1971}.

\begin{Proposition}[Spectrum of Sturmian numbers]
For distinct $\alpha_-,\alpha_+\in \R\P^1$, the following are equivalent:
\begin{itemize}[noitemsep]
    \item The geodesic $(\alpha_-,\alpha_+) \subset \M'$ is simple and closed.
    \item $\alpha_\pm$ are conjugate Markov quadratics.
    \item We have $\Cohn(\alpha_-,\alpha_+)>4$.
\end{itemize} 
and the following are equivalent:
\begin{itemize}[noitemsep]
    \item The geodesic $(\alpha_-,\alpha_+) \subset \M'$ is simple but not closed.
    \item $\alpha_\pm$ are conjugate Sturmian numbers which are non-Markov.
    \item We have $\Cohn(\alpha_-,\alpha_+) = 4$.
\end{itemize} 
Moreover, when $(\alpha_-,\alpha_+)\subset \M'$ simple, we have $\Cohn(\alpha_-,\alpha_+)=4\coth\left(\tfrac{1}{2}\ell_{\M'}(\alpha) \right)$.
\end{Proposition}

\begin{Question}[simple in $\M''$]
    What are the Lagrange or Markov-Cohn spectra for the endpoints of lifts $(\alpha_-,\alpha_+) \subset \HP$ of simple or simple closed geodesics of $\M''$?
\end{Question}

\subsection{Transcendence of the Inshection and hexponential}

We define the subset $\Markov\subset \Sturm$ of \emph{Markov} numbers as those arising from rational slopes $(x,y) \in \S H_1(\M';\Z)$, namely those with periodic $X^{\pm 1}\& Y^{\pm 1}$ sequences.
These coincide with the Markov quadratic irrationals \cite{Series_Geo-Markov-Num_1985}, and the involution $\sigma$ restricts to the action by Galois conjugation.

\begin{Theorem}[Transcendence of the Inshection map]
\label{thm:transcendence-hexp}
    For $(x,y) \in \S H_1(\M;\R)$, 
    \begin{itemize}[noitemsep]
        \item if $(x,y) \in \S H_1(\M';\Z)$ is rational then $\InSh(x,y)\subset \Markov$ are quadratic
        \item if $(x,y) \notin \S H_1(\M';\Z)$ is irrational then $\InSh(x,y)\subset \Sturm \setminus \Markov$ are transcendental
    \end{itemize}
\end{Theorem}

\begin{proof}
If $(x,y) \in \S H_1(\M';\Z)$ is rational then every number in $\InSh(x,y)$ has a continued fraction expansion which is eventually periodic, so it is quadratic.
If $(x,y) \notin \S H_1(\M';\Z)$ is irrational, then every number in $\InSh(x,y)$ has a continued fraction expansion which is an eventually aperiodic Sturmian sequence over $\{1,2\}$, so it is transcendental by \cite[Proposition 3 and Theorem 7]{ADQZ_transcendence-Sturmian_2001}.

The strategy in \cite{ADQZ_transcendence-Sturmian_2001} consists of showing that a non-Markov Sturmian number is very well approximated by Markov numbers, well enough to satisfy the hypotheses of a theorem of Schmidt \cite{Schmidt_simult-approxim-algebraic-by-rational_1967}, which can be stated as follows. The height of a quadratic irrational $\xi$ with minimal polynomial $ax^2+bx+c\in \Z[x]$ is defined by $H(\xi)=\gcd(\lvert a \rvert ,\lvert b\rvert ,\lvert c\rvert )$.
A number $s\in \R\setminus \Q$ which is not quadratic must be transcendental, if for some number $\beta>3$ there exists infinitely many quadratic irrationals $\xi_k$ such that $\lvert \xi-\xi_k\rvert <H(\xi_k)^{-\beta}$.
\end{proof}

\begin{Question}
\sloppy To study the diophantine properties of the ray compactification $\Shexp \colon \Radial \to \S H_1(\M',\R)$ and its inverse $\InSh \colon \S H_1(\M',\R) \to \Sturm$, one may try to find continued fractions for $\tan \circ \Shexp \colon \Radial \to \R\R^1$ and $\InSh \circ \arctan \colon \R\P^1 \to \Sturm$.
\end{Question}


We thank Bill Duke for communicating to us the following continued fraction expansion of $\hexp$ in terms of the modular function $\lambda \colon \HP \to \C\setminus \{0,1,\infty\}$ (uniformizing the congruence cover $\M(2)$ of $\M$ associated to the congruence subgroup $\Gamma(2)$ of $\Gamma$).

\begin{Proposition}[Continued fraction for $\hexp(\lambda)$]
In terms of $\lambda = \lambda(\tau)$ we have:
\begin{equation*}
    \int_\infty^\tau \eta^4(z)dz
    =\tfrac{3}{i2\pi}\left(\tfrac{1}{2}\lambda(1-\lambda)\right)^{1/3}\cdot \frac{1}{1-\frac{n_1\lambda}{1-\frac{n_2 \lambda}{1 -\dots}}}
\end{equation*}
where for $k\in \N$ we have  $n_{2k+1}=\frac{(k+1/3)}{2(2k+1/3)}$ and for $k\in \N^*$ we have $n_{2k}=\frac{k}{2(2k+1/3)}$.
\end{Proposition}

\begin{proof}
We first recall the hypergeometric function ${}_2F_1(a,b,c;z)$ defined for $a,b,c\in \C$ with $-c\notin \N$ by the following power series expansion which converges for $\lvert z\rvert < 1$:
\begin{equation*}
    {}_2F_1(a,b,c;z) = \sum_{n=0}^\infty \frac{(a;n)(b;n)}{(c;n)(1;n)} z^n 
    \qquad \mathrm{where} \quad (k;n)=k\dots (k+n-1).
\end{equation*}
It satisfies $\left[(a+D_z)(b+D_z)-(c+D_z)(1+D_z)z^{-1} \right] F(z)=0$ where $D_z = z\frac{d}{dz}$, that is the differential equation of order $2$ with regular singularities $\{0,1,\infty\}$:
\begin{equation*}
    z(1-z) \tfrac{d^2}{dz^2} F+(c-(a+b+1)z) \tfrac{d}{dz} F(z) - ab F(z) = 0
\end{equation*}
Moreover there is the Gauss continued fraction discovered by Gauss, and whose convergence for $\lvert z \rvert < 1$ was settled in \cite{Van-Vleck_convergence-Gauss-fracont_1901} and has been revisited in \cite{Ifantis-Panagopoulos_convergence-fracont-revised_2001}:
\begin{equation*}
    \frac{{}_2F_1(a,b+1,c+1;z)}{{}_2F_1(a,b,c;z)} = \frac{1}{1-\frac{n_1 z}{1-\frac{n_2 z}{1 -\dots}}}
\end{equation*}
where $n_{2k+1}=\frac{(a-c-k)(b+k)}{(c+2k)(c+2k+1)}$ for $k\in \N$ and $n_{2k}=\frac{(b-c-k)(a+k)}{(c+2k-1)(c+2k)}$ for $k\in \N^*$.

Using the computations from \cite[pages 434 to 438]{Kleban-Zagier_Crossing-probabilities-modular-forms_2002} we find the expression 
\begin{equation*}
    \int_\infty^\tau \eta^4(z)dz
    =\tfrac{3}{i2\pi}\left(\tfrac{1}{2}\lambda\right)^{1/3}\cdot {}_2F_1(1/3,2/3,4/3;\lambda)
\end{equation*}
and dividing by $(1-\lambda)^{1/3}={}_2F_1(1/3,-1/3,1/3;\lambda)$ yields the continued fraction.
\end{proof}

\begin{Remark}[Convergence]
The previous parametrisation may appear surprising, as the symmetry group $\Gamma''$ of $\hexp$ is smaller than the symmetry group $\Gamma(2)$ of $\lambda$.
To be more precise $\Gamma''$ is commensurable to $\Gamma(2)'$, and we have $\Gamma(2)/\Gamma(2)'=\Z^2$.

In fact the continued fraction expression for $\hexp(\lambda)$ converges only for $\lvert \lambda \rvert < 1$.
However, it admits an analytic continuation around its three singularities $\{0,1,\infty\}$ and we deduce that $\Gamma(2)/(\Gamma(2)\cap \Gamma'')= \Z^2\rtimes \Z/2$ must coincide with the monodromy groug of $\hexp(\lambda)=C\tfrac{3}{i2\pi}\left(\lambda/2\right)^{1/3}\cdot{}_2F_1(1/3,2/3,4/3;\lambda)$.
\end{Remark}

\begin{Remark}[Monodromy]
    Note that the function $\lambda$ is an algebraic function of $\J$, since the cover $\HP/\Gamma(2) \to \HP/\Gamma$ has solvable Galois group $\mathfrak{S}_3$.
    We also know that that $\hexp$ and $\lambda$ are algebraically independant, since the cover $\HP/\Gamma'' \to \HP/\Gamma$ has infinite Galois group $\Z^2\rtimes \Z/6$.
    However this group is solvable.

    This shows that the expression of $\hexp$ in terms of $\lambda$ is in some sense solvable with a Galois group essentially equal to $\Z^2$.
    It would be interesting to relate the solvability of this quadratic Galois group with the quadratic growth in the complexity of the coefficients $n_k$ in the continued fraction expansion for $\hexp(\lambda)$.
\end{Remark}

\appendix

\section{Simple geodesics in other Riemann surfaces}

\label{sec:other-surfaces}

\subsection{Arithmetical and dynamical simplicity of geodesics}

Let us announce a conjecture, part of which we will prove in a forthcoming paper.

\begin{Conjecture}[Arithmetic-Dynamic simplicity]
\label{thm:low-complexity-dichotomy}
Consider a lattice $F\subset \PSL_2(\R)$ with invariant trace field $\K=\Q(\{\disc(\gamma)\mid \gamma \in F\})$ and assume that the entries of $F$ belong to $\K$ (this may be achieved by a conjugacy in $\PSL_2(\R)$). Let $S=F\backslash \HP$.

Fix a complete simple geodesic $\alpha \subset S$ and choose any lift $(\alpha_-,\alpha_+)\subset \HP$ where $\alpha_\pm \in \R\P^1$ are its endpoint.
We conjecture that the following trichotomy holds:
\begin{itemize}[noitemsep, align=right]
    \item[$1$:] If $\alpha$ is forward asymptotic to a cusp, then $\alpha_+$ belongs to $\K$.
    \item[$2$:] If $\alpha$ is forward asymptotic to a closed geodesic, then $\alpha_+$ is quadratic over $\K$.
    \item[$\infty$:] Otherwise, $\alpha_+$ is transcendental over $\K$.
\end{itemize}
Of course, this would imply the same result exchanging $\alpha_+$ for $\alpha_-$ and forwards asymptotic with backwards asymptotic.

Note that the statement does not depend on the choice of the lift, since all the lifts differ by elements in $F\subset \PSL_2(\K)$.

The points $1$ and $2$ are easy to prove, the heart of the conjecture lies in $\infty$.
\end{Conjecture}

\begin{Remark}[Heuristics on arithmetic-dynamic complexity]
This theorem and the ideas involved in the proof lead to the following heuristic principle.
\begin{quote}
    Objects with low dynamical and arithmetic complexity follow a dichotomy: they present either very algebraic features or very transcendent features.
\end{quote}


Such dichotomy results appear as a leitmotiv in transcendental function theory.
One expects that a power series with bounded dynamical complexity (say on the dimension of the vector spaces generated by its derivatives) and bounded arithmetic complexity (say on the growth of the algebraic heights of its coefficients) either defines an algebraic function and its values at algebraic numbers are all algebraic or else defines a very transcendent function in the sense that its values at almost all algebraic numbers are transcendent.

Striking examples are the Siegel-Shidlovskii theorem improved by Beukers in \cite{Beukers_Refined-Siegel-Shidlovskii-Theorem_2004} or the Shneider-Lang theorem which is presented in  \cite[(3.1)]{Lang_Transcendental-numbers-diophantine-approximations_1971}.

This leitmotiv will appear in the subsections.
\end{Remark}

\subsection{Real algebraic cycles in complex algebraic curves}

A theorem of Schneider \cite{Schneider_transcendant-functions_1949} about the modular function $\J\colon \HP \to \C$ says that the value $\J(\tau)$ is algebraic (namely the elliptic curve $E_\tau$ is defined over $\overline{\Q}$) if and only if its argument $\tau$ is complex quadratic (namely the elliptic curve $E_\tau$ has a period-lattice admitting complex multiplication).

Our cusp compactification $\partial \hexp$ and ray-compactification $\Shexp$ give some sense to the limits of the modular function $\J$ at the boundary.
One may thus view our theorem \ref{thm:transcendence-hexp} as providing the algebraic-transcendence dichotomy for the values of $\InSh$ (the inverse to the ray-compactification $\Shexp$ of $\hexp$).

The theorem of Schneider has been generalised to the Siegel moduli space of principally polarized abelian varieties, parametrized as the space of complex tori $E_\tau = \C^{g}/\Z^g+\tau \Z^g$ where $\tau$ is now a symmetric matrix whose imaginary part is positive definite.
Note that a complex algebraic curve $S$ of genus $g$ defined over $\overline{\Q}$ has a Jacobian variety $H_1(S;\C)/H_1(S;\Z)$ of dimension $2g$ which defines an algebraic point in that moduli space.

Instead of stating this generalisation, let us mention the analytic subgroup theorem of \cite{Wustholz_analytic-subgroup-theorem_1989} on which it relies. Let $G$ and $G'$ be commutative algebraic groups defined over the field $\overline{\Q}$ of algebraic numbers, and of positive dimensions. Let $\varphi \colon G'(\C) \to G(\C)$ be an analytic homomorphism, which is defined over $\overline{\Q}$ (in the sense that the induced map $d\varphi \colon TG'(\C) \to TG(\C)$ on tangent spaces is compatible with the $\overline{\Q}$ structures). If $\varphi(G')(\overline{\Q})\ne \emptyset$ then $\varphi(G')$ contains an algebraic subgroup of $G$ which is defined over $\overline{\Q}$ and of positive dimension.

Consider a lattice $F\subset \PSL_2(\R)$, and the complex algebraic curve $S=F\backslash \HP$. It is known that the lattice $F$ is arithmetic (in the sense of Margulis or Shimura \cite{Mochizuki_Correspondences-hyperbolic-curves_1998}) if and only if the curve $S$ (or its Jacobian variety $\operatorname{Jac(S)}$) is defined over $\overline{\Q}$ with periods in $\overline{\Q}$, if and only if the the curve $S$ (or $\operatorname{Jac}(S)$) has complex multiplication. 

\begin{Conjecture}
Consider a lattice $F\subset \PSL_2(\R)$ such that the complex algebraic curve $S=F\backslash \HP$ (or equivalently its Jacobian) is defined over $\overline{\Q}$.

Consider a simple geodesic $\alpha\subset S$ and denote its lift by $\Tilde{\alpha}\subset \HP$ which is an $F$-invariant set of geodesics.
We conjecture that the following are equivalent:
\begin{itemize}[noitemsep]
    \item the simple geodesic $\alpha$ is closed
    \item the simple geodesic $\alpha$ contains an algebraic point of $S$
    \item the endpoints $\partial\Tilde{\alpha}\subset \R\P^1$ are algebraic
\end{itemize}
Note that for simple closed curves $\alpha\subset S$ the lifts have endpoint $\partial \Tilde{\alpha}\subset \R\P^1$ which are quadratic over the invariant trace field $\K = \Q\{\disc(A)\mid A \in F\}$ of $F$.

\end{Conjecture}

\subsection{Complex multiplication and real multiplication}

Now suppose that $\tau\in \HP$ is complex quadratic, namely that $\J(\tau)$ is algebraic.

The main theorem of class field theory for complex quadratic fields states that for such pairs $(\tau, \J(\tau)) \in \overline{\Q}$, the quadratic field $\Q(\tau)$ has maximal unramified abelian extension generated by $\J(\tau)$, whose maximal abelian extension is further generated by the special values of certain $\Lambda_\tau$-elliptic functions at torsion points in $\C/\Lambda_\tau$.

Observe that $\J\colon \HP \to \M$ factors through $\hexp \colon \HP \to \M''$, and the aforementionned arithmetic results about $\J(\tau)$ also hold for $\hexp(\tau)$.

It was Kronecker's Jugendtraum (dream of youth) to find a generalisation of  for real quadratic fields, which eventually became Hilbert's 12-th problem. 
One may consult \cite{Schappacher_history-Hilbert-12_1998} for the history of that still unsolved problem.
Let us comment on some approaches which have been proposed, and relate them to this work.

Consider a quadratic irrational $\alpha_+ \in \R$ with Galois conjugate $\alpha_-$ and let $A\in \Gamma$ be the primitive matrix with axis $(\alpha_-,\alpha_+)\subset \HP$.

The ray compactification $\Shexp \colon \Radial \to \R/2\pi\Z$ will not suffice to understand the Hilbert class field of $\alpha$.
Indeed $\alpha \in \Radial \iff A^6\notin \Gamma''$, and $\alpha\in \Markov \implies \Shexp(\alpha) \in \Q$.

However the geodesic $\hexp(\alpha',\alpha)\subset \C\setminus \Lambda $ contains all the information about $\alpha$ and we can hope to generate the class field of $\Q(\alpha)$ by integrating elliptic functions in $\Q[\wp_\Lambda, \wp'_\Lambda]$ along the geodesic or by summing quantities measuring how close the geodesic approaches the lattice point $\partial\hexp(\Q\P^1) \subset \M''$.
These approaches relate to the ideas in \cite{Choie-Zagier_Rational-period-functions-PSL2Z_1993, Duke-Imamoglu-Toth_cycle-integrals-j-function_2011} and \cite{Bernard-Gendron_modular-invariant-quantum-tori_2014} respectively.

\begin{Remark}
It would be interesting to relate the periods of $\hexp$ and the values of $\Shexp$ to the periods and special values of certain modular cocycles introduced in \cite[Section 2.4 and 2.5]{Manin_Real-multiplication-noncommutative-geometry_2004} and developped in.

More precisely, given $A\in \Gamma'$ with fixed points $\alpha\pm \in \R\P^1$, can we relate the period $\int_z^{Az} \hexp(\tau) d\tau$ for $z\in(\alpha_-, \alpha_+)$, or the value $\Shexp(\alpha_+)-\Shexp(\alpha_-)$ to the value of the modular cocycles $\mathfrak{w}_A(\alpha_+)$ associated to constructed by Knopp.
\end{Remark}



\newpage

\section{On the Fourier coefficients of the hexponential}

\label{sec:more-on-psi}

Let us collect some well known facts about the sequence $\psi(n)$.
Our main reference are \cite{Ligozat_Courbes-modulaires-genre-1_1975, Serre_lacunarite-puissances-eta_1985} and \cite{lmfdb}.

\subsection{Modularity and growth of the coefficients}

The $\psi(n)$ correspond to the Fourier coefficients of the weight $2$ modular form $\eta^4(6\tau)$, which is \cite[\href{https://www.lmfdb.org/ModularForm/GL2/Q/holomorphic/36/2/a/a}{Modular Form 32.2.a.a}]{lmfdb}, the unique normalised Hecke cusp form for the group $\Gamma_0(36)= \{\begin{psmallmatrix}
a&b\\c&d \end{psmallmatrix}\in \Gamma \mid c\equiv 0 \bmod{36}\}$.
It has an associated Dirichlet series $L(\psi,s)=\sum_1^\infty \tfrac{1}{n^s}\psi(n)$. 
Let us recall the following from \cite{Ligozat_Courbes-modulaires-genre-1_1975}.

\begin{enumerate}[align=left]
    \item[Théorème B (2.2.3)] The Dirichlet series $L(\psi,s)$ attached to the group $\Gamma_0(36)$ coincides with the $L$-series attached to the elliptic modular curve $X_0(36)$ (with Weierstrass form $y^2=x^3-1$, having complex-multiplication by $\Z[j]$).
    
    \item[Théorème C (2.3.2)] The function $\Xi(s)=\left(\tfrac{6}{2\pi}\right)^s \Upgamma(s) L(\psi,s)$ admits an analytic continuation to an entire function which is bounded in the vertical strip and satisfies the functional equation $\Xi(s)=\Xi(2-s)$. 

    
    \item[Théorème F (6.3.3)] The analytic continuation of the function $L(\psi,s)$ converges at $s=1$ to 
    a constant divided by the square of the order of its group of rational points $X_0(36)(\Q)$, which is finite by Théorème E.
\end{enumerate}

In particular the sequence $\psi(n)$ satisfies the Ramanujan-Petersson conjecture \cite{Deligne_SB-Formes-modulaires-representations-l-adiques_1971} saying that $\lvert \psi(p) \rvert \le 2\sqrt{p}$ for all prime $p\in \N$.

\subsection{L-series and Euler product}

The ring of integers $\Z[j]$ of the quadratic field $\Q(\sqrt{3})$ has Euclidean algorithms, in particular it is a principal ideal domain, hence a unique factorisation domain.

The ramification pattern of a prime $p\in \Z$ in $\Z[j]$ is given by its residue $p\bmod{3}$ as follows.
The prime $3$ is ramified as $3=-(1+j)^2$.
When $p=-1\bmod{3}$ it is inert, namely it remains prime in $\Z[j]$, and has norm $p^2$.
When $p=+1\bmod{3}$ it splits in $\Z[j]$ as a product $p=\pp \Bar{\pp}$ of conjugate primes with $\Norm(\pp)=p$.


The \emph{Hecke Grossencharacter} $\ch$ of discriminant $-3$ and conductor $2\sqrt{-3}$ is a character on the multiplicative group of fractional ideals of $\Q(j)$ coprime to the conductor $2\sqrt{-3}$ defined as follows. 
For an ideal $\mathfrak{a}\subset \Z[j]$ coprime to $2\sqrt{-3}$, we have $\ch(\mathfrak{a})=\alpha$ where $\alpha$ is the unique generator of $\mathfrak{a}$ such that $\alpha = 1 \bmod{2\sqrt{-3}}$.
In particular for a prime $p\in \Z$ coprime to $6$ we have $\ch(p)=\pm 1$ according to whether $p=\pm 1\bmod{3}$.

We have the $L$-series representation where the product is over all primes $\pp\in \Z[j]$ coprime to the conductor $2\sqrt{-3}$:
\begin{equation*}
L(\psi,s)=\sum_1^\infty \frac{\psi(n)}{n^s} 
= \prod_\pp \frac{1}{1-\ch(p)\Norm(p)^{s}}
\end{equation*}
One may rewrite this as the product over all primes $p\in \Z$ coprime to $6$:
\begin{equation*} 
L(\psi,s)=\sum_1^\infty \frac{\psi(n)}{n^s} = \prod_{p} L_p(\psi,s)
\end{equation*}
where the local factors are given by $L_p(\psi,s)=(1+p\cdot p^{-2s})^{-1}$ if $p\ne 2$ is inert and $L_p(\psi,s)=(1-(\ch(\pp)+\ch(\Bar{\pp})) p^{-s}+p\cdot p^{-2s})^{-1}$ if $p=\pp \Bar{\pp}$ is split.
Hence the multiplicative function $\psi \colon \N \to \Z$ is determined on powers of primes $p^n$ for $n\in \N$ by:
\begin{enumerate}[noitemsep]
    \item[-] If $p\in \{2,3\}$ then $\psi(p^n)=0$. 
    \item[-] If $p\ne 2$ is inert then $\psi(p^n)\in \{0,1\}$ is equal to the parity $n\bmod{2}$.
    \item[-] If $p=\pp\Bar{\pp}$ splits then $\psi(p^n)=\ch(\pp)^n+ \ch(\pp)^{n-1} \ch(\Bar{\pp}) +\dots+ \ch(\pp) \ch(\Bar{\pp})^{n-1} +\ch(\Bar{\pp})^n$
\end{enumerate}

\subsection{Lacunarity}

The \cite[Section 2.2]{Serre_lacunarite-puissances-eta_1985} characterises the values $n\in \N$ such that $\psi(n)= 0$ in terms of the ramification of $n$ in the ring of integers $\Z[j]$ of the quadratic field $\Q(j)$: 
\begin{quote}
For $n\in \N^*$ we have $\psi(n)=0$ if and only if there exists a prime $\mathfrak{p}\in \Z[j]$ with $\mathfrak{p}\equiv -1 \bmod{3}$ such that the $\mathfrak{p}$-adic valuation of $n$ is odd.
\end{quote}

Then Serre deduces \cite[Théorème 2(ii)]{Serre_lacunarite-puissances-eta_1985} from an application of the Chebotarev density theorem in \cite[Proposition 18]{Serre_applications-densite-Chebotarev_1981}, that there exists $c_4\in \R_+^*$ (for which an explicit value is given) such that the following asymptotic holds: 
\begin{equation*}
    \tfrac{1}{N} \Card\{n\in \N\mid \psi(n)\ne 0, \, n\le N\} \sim \tfrac{c_4}{\log(N)^{1/2}}
\end{equation*}
In particular the sequence $\psi$ is lacunary in the sense that the left hand side converges to $0$, which was known by \cite{Ramanujan_certain-arithmetic-functions_1916}.
As Remarked in \cite[page 217]{Serre_lacunarite-puissances-eta_1985}, this lacunarity is barely pronounced, which partly explains the following.

\begin{Remark}[Experimental]
\label{rem:experimental-slow-convergence}
We computed the first $20000$ terms of the partial series $\partial \hexp(r)$ for $r\in \Q$ using \href{http://pari.math.u-bordeaux.fr/}{PARI GP}, and observed that they start to coalesce in some region, but then oscillate in that small region, at an amplitude which seemed to be neither approaching $0$ (although the proof reveals they must be), nor monotonous. 
\end{Remark}

\section*{Acknowledgements}

We wish to thank Rich Schwartz, Peter Doyle and John Voight for the encouragement, discussions and relevant remarks, as well as Winnie Li and Fang-Ting Tu for their continuous interest and effort in reading and replying to our emails. We are also indebted to Keith Conrad for answering our question \cite{CLS_Quest-MathOverflow_2024}. 

This material is based upon work supported by the National Science Foundation under Grant No. DMS-2247553.

\bibliographystyle{alpha} 
\bibliography{paratext/ref_Hexp.bib}

\end{document}